\documentclass[12pt]{article}
\usepackage{amsmath}
\usepackage{amsthm}
\usepackage{amsbsy}
\usepackage{amssymb}
\usepackage{amscd}
\usepackage{a4wide}
\usepackage{vmargin}
\usepackage{epic}

\newtheorem{lemma}{Lemma}
\newtheorem{theorem}{Theorem}
\newtheorem{statement}{Proposition}
\newtheorem{definition}{Definition}

\newcommand{\set}[1]{\left\{#1\right\}}

\DeclareMathOperator{\rad}{Rad} 
 \DeclareMathOperator{\rep}{Rep}
 
 \DeclareMathOperator{\ob}{Ob}
 \DeclareMathOperator{\diag}{diag}
\DeclareMathOperator{\rk}{rk}

\begin{document}

\title{ The Spectral  Problem and Algebras
Associated with Extended Dynkin Graphs.}

\author{Stanislav Popovych \footnote{
E-mail address: stanislp\symbol{64}math.chalmers.se. Tel:+46 0768701227 Fax: +46 31 161973} \\ \vspace{0.1cm}\\
{\footnotesize Department of Mathematics, Chalmers University of
Technology,
SE-412 96 G\"oteborg, Sweden} \\
{\footnotesize \sl
stanislp\symbol{64}math.chalmers.se}\\ \vspace{0.3cm}\\
  }


\date{}
\maketitle{}

\numberwithin{equation}{section}
 \abstract{  The  Spectral Problem is to describe
  possible spectra $\sigma (A_j)$ for an irreducible $n$-tuple
   of Hermitian operators s.t. $A_1+\ldots+A_n$ is a scalar operator.
   In case when $m_j= | \sigma (A_j)|$ are finite and a rooted tree
    ${\rm T}_{m_1,\ldots, m_n}$ with $n$ branches of lengths
    $m_1, \ldots, m_n$ is a Dynkin graph the explicit
    answer to the Spectral Problem was given recently
    by S.~A. Kruglyak, S.~V. Popovych, and Yu.~S. Samo\v\i{}lenko in~\cite{KPS1}. In present work the solution of  the Spectral
    Problem for all star-shaped simply laced extended Dynkin graphs, i.e.
    when  $(m_1, \ldots, m_n) \in \{ (2,2,2,2), (3,3,3), (4,4,2)$, $(6,3,2)\}$ is presented.

 }
\vspace{10mm}

\newpage

\noindent {\Large \bf Introduction.} \vspace{5mm}

 {
                         Let  $A_1$, $A_2$, $\ldots$, $A_n$  be Hermitian  $m\times m$
 matrices with given
eigenvalues: $\tau(A_j)= \{\lambda_1(A_j)\ge \lambda_2(A_j)\ge
\ldots \ge \lambda_m(A_j) \}$. The well-known classical problem
about the spectrum of a sum of two Hermitian matrices (Horn's
problem) is to describe possible values of
 $\tau(A_1), \tau(A_2), \tau(A_3)$  such that $A_1+A_2=A_3$. In
 more symmetric setting one can seek for a connection between $\tau(A_1), \tau(A_2), \ldots,
 \tau(A_n)$ necessary and sufficient for the existence of Hermitian
 operators such that $A_1+A_2+\ldots+A_n=\gamma I$ for a fixed $\gamma\in
 \mathbb{R}$.

A recent solution of this problem  (see ~\cite{fulton, klychko}
 and others) gives a complete description of
possible $\tau(A_1), \tau(A_2), \ldots,
 \tau(A_n)$ in terms of linear
inequalities conjectured by A. Horn.

 A modification of Horn's problem called henceforth the {\it spectral problem}  was considered in \cite{KPS1, KPS}.
Let $A_1$, $A_2$, $\ldots$, $A_n$ be bounded linear Hermitian operators  on a separable Hilbert space $H$. For an operator
$X$ denote by $\sigma(X)$ its spectrum. Given  $M_1, M_2, \ldots,
M_n$  closed finite subsets of $\mathbb{R}$ and
$\gamma\in\mathbb{R}$ the problem is
to determine whether there are Hermitian operators $A_1$, $A_2$,
$\ldots$, $A_n$ on $H$ such that $\sigma(A_j)\subseteq M_j$
 ($1\le j\le n$), and
$A_1+A_2+\ldots+A_n=\gamma I$.

The essential difference with Horn's classical  problem is that we
do not fix the dimension of  $H$   (it may be finite or infinite)
  and the spectral multiplicities.

 The  Spectral Problem can be stated  in terms of
$*$-representations of $*$-algebras.  Namely, let
$\alpha^{(j)}=(\alpha_1^{(j)},\alpha_2^{(j)},\ldots,
\alpha_{m_j}^{(j)})$ ($1\le j\le n$)  be vectors with positive
strictly decreasing coefficients.  Put  $M_j=\alpha^{(j)}\cup \{0\}$. Let us
consider the associative algebra defined by the following
generators and relations \begin{gather*} \mathcal{A}_{M_1, \ldots, M_n, \gamma} = \langle A_1, \ldots, A_n |
A_1+ \ldots +A_n = \gamma e,\\
 P_1(A_1) = 0, \ldots, P_n(A_n)=0 \rangle
\end{gather*}
where $P_k$ is a polynomial with simple roots from the set $M_k$.
 This is a  $*$-algebra, if we declare all generators $A_j$ to be
self-adjoint. A $*$-representation $\pi$ of $\mathcal{A}_{M_1,  \ldots, M_n,
\gamma}$  determines  an $n$-tuple  of non-negative operators
$A^{(j)}= \pi(A_j)$ with spectral decomposition
$\sum_{k=1}^{m_i} \alpha_k^{(i)} P_k^{(i)}$. Let
$P_0^{(i)} = I- \sum_{k=1}^{m_i} P_k^{(i)}$ be the spectral projection corresponding to zero eigenvalue.
 Clearly  any
$n$-tuple of Hermitian matrices $A^{(1)}$, $A^{(2)},\ldots$,
$A^{(n)}$ such that $A^{(1)}+\ldots+ A^{(n)}=\gamma I$ and
$\sigma(A^{(j)})\subseteq M_j$ determines a representation $\pi$
of $\mathcal{A}_{M_1,  \ldots, M_n, \gamma}$.
 So in terms of $*$-representations, the
spectral problem is a problem to describe  the set $\Sigma_{m_1,m_2,\ldots,m_n}$
of the parameters $\alpha_k^{(j)}$, $\gamma$ for which there exist
$*$-representations of $\mathcal{A}_{M_1,  \ldots, M_n, \gamma}$.
  We will call a  $*$-representation
$\pi$  of the algebra
 $\mathcal{A}_{M_1,  \ldots, M_n, \gamma}$  {\it non-degenerate}
if  $\pi(p_k^{(j)})\not=0$ for all $0\le k \le m_j$ and $1\le j \le n$.

  With an  integer vector    $(m_1,\ldots,m_n)$ we
will
  associate a non-oriented  star-shape  graph  $G$
  with $n$ branches of the lengths   $m_1$, $m_2$, $\ldots$,
  $m_n$  stemming from a  single root. The graph $G$ and vector
  $\chi= (\alpha_1^{(1)}, \alpha_2^{(1)}$, $\ldots$, $\alpha_{m_1}^{(1)},
  \alpha_{1}^{(2)}$, $\ldots, \alpha_{m_2}^{(2)}$, $\ldots,
  \alpha_{1}^{(n)}$, $\ldots$, $\alpha_{m_n}^{(n)}$, $\gamma)$ completely
  determine the algebra $\mathcal{A}_{M_1,  \ldots, M_n,
  \gamma}$   so we will write $\mathcal{A}_{G,\chi}$ instead of $\mathcal{A}_{M_1,  \ldots, M_n,
  \gamma}$.  Henceforth we will denote the set $\Sigma_{m_1,\ldots,m_n}$
  by  $\Sigma(G)$ where $G$ is the tree mentioned above.

By  $\overline{\rep} \mathcal{A}_{G, \chi}$
 we will denote the full subcategory of non-degenerate
representations of the category ${\rep} \mathcal{A}_{G, \chi}$ of
$*$-representations of the $*$-algebra $\mathcal{A}_{G, \chi}$ on
 Hilbert spaces.

To solve the Spectral Problem it is enough to consider only non-degenerate representations. 
 Consider the following set:
    $\Sigma_{m_1,\ldots,m_n}^{n.-d.}=
  \{ ((\alpha_k^{(j)})_{jk},\gamma)|$ { \it there is a non-degenerate $*$-representation of }
  $\mathcal{A}_{M_1,  \ldots, M_n, \gamma}\}$;
  which depends only on   $(m_1$, $\ldots$, $m_n)$.
 Every irreducible representation of the algebra $\mathcal{A}_{M_1,
M_2, \ldots, M_n, \gamma}$ is an irreducible non-degenerate
*-representation of an algebra $\mathcal{A}_{\widetilde{M}_1,
\ldots, \widetilde{M}_n, \gamma}$ for some subsets
$\widetilde{M}_j\subset M_j$. Hence $(M_1,  \ldots, M_n,
\gamma)\in \Sigma_{m_1,\ldots, m_n}$ if there exist
$(\widetilde{M}_1, \ldots, \widetilde{M}_n,\gamma)\in
\Sigma_{|\widetilde{M}_1|,\ldots, |\widetilde{M}_n|}^{n.-d.}$.
Thus the description of $\Sigma_{m_1,\ldots, m_n}$ follows from
the description of $\Sigma_{k_1,\ldots, k_n}^{n.-d.}$ where
$k_j\le m_j, 1\le j\le n$.

If the  graph $G$ is a Dynkin graph or extended Dynkin graph the
problem is greatly simplified. The algebras $\mathcal{A}_{G,\chi}$ associated with Dynkin graphs  (resp.
extended Dynkin graphs) have a more simple structure than in other
cases. In particular, the algebras $\mathcal{A}_{G,\chi}$ are finite dimensional (resp. have polynomial
growth) if and only if the associated graph is a Dynkin graph
(resp. an extended Dynkin graph) (see~\cite{melit}).

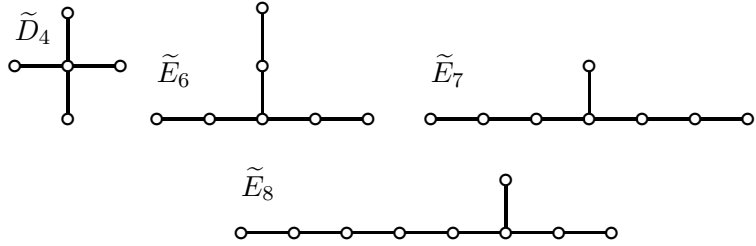
\setlength{\unitlength}{2pt}
\begin{center}
\begin{picture}(25,25)(-1,-1)
\thicklines
\put(-1,15){$\widetilde{D}_4$}
 \drawline(9,11)(9,19)
 \drawline(9,9)(9,1)\put(9,0){\circle{2}}
 \multiputlist(0,10)(10,0){\circle{2},\circle{2},\circle{2}}
 \multiputlist(10,20)(0,10){\circle{2}}
 \drawline(0,10)(8,10)
\drawline(10,10)(18,10)
\end{picture}
\begin{picture}(50,25)(-1,-1)
\thicklines
\put(-1,7){$\widetilde{E}_6$}
 \drawline(19,1)(19,9)
 \drawline(19,11)(19,20)\put(19,21){\circle{2}}
\multiputlist(0,0)(10,0){\circle{2},\circle{2},\circle{2},\circle{2},\circle{2}}
 \multiputlist(20,10)(0,10){\circle{2}}
 \drawline(0,0)(8,0)
\drawline(10,0)(18,0)
\drawline(20,0)(28,0)\drawline(30,0)(38,0)
\end{picture}
\begin{picture}(69,15)(-1,-1)
\thicklines
\put(-1,7){$\widetilde{E}_7$}
 \drawline(29,1)(29,9)
\multiputlist(0,0)(10,0){\circle{2},\circle{2},\circle{2},\circle{2},\circle{2},\circle{2},\circle{2}}
 \put(29,10){\circle{2}}
 \drawline(0,0)(8,0)
\drawline(10,0)(18,0)
\drawline(20,0)(28,0)\drawline(30,0)(38,0)\drawline(40,0)(48,0)\drawline(50,0)(58,0)
\end{picture}
\end{center}
\begin{center}
\begin{picture}(62,15)(-1,-1)
\thicklines
\put(-1,7){$\widetilde{E}_8$}
 \drawline(49,1)(49,9)
\multiputlist(0,0)(10,0){\circle{2},\circle{2},\circle{2},\circle{2},\circle{2},\circle{2},\circle{2},\circle{2}}
 \put(49,10){\circle{2}}
 \drawline(0,0)(8,0)
\drawline(10,0)(18,0)
\drawline(20,0)(28,0)\drawline(30,0)(38,0)\drawline(40,0)(48,0)\drawline(50,0)(58,0)\drawline(60,0)(68,0)
\end{picture}
\end{center}

It is known (see \cite{ostrovskij1}, \cite{melit1}) that all irreducible $*$-representations of an algebra $\mathcal{A}_{G, \chi}$ with Dynkin or extended Dynkin graph  $G$ are finite-dimensional. Thus for such graphs it is enough to consider Spectral Problem in finite-dimensional Hilbert spaces. 
 
It was shown in~\cite[Theorem~3]{ostrovskij} that if $G$ is not a Dynkin nor an extended Dynkin graph then algebra $\mathcal{A}_{G, \chi}$ has an infinite-dimensional irreducible representation for some value of $\chi$. 

For algebras associated with  Dynkin graph $G$ an explicit solution
of the Spectral Problem was obtained in~\cite{KPS1}. In the present
paper we present an explicit solutions of the Spectral Problem for
all star-shaped  simply-laced extended Dynkin graphs (see picture
above).

Let us explain the solution.   Let  $G$  be a  star-shaped extended
simply-laced Dynkin graph. For simplicity we consider below only
graphs with three branches. With each irreducible solution of the
Spectral Problem associated with graph $G$ we have the corresponding family of
self-adjoint operators $A^{(j)}$ ($j=1,2,3$) on a Hilbert space $H$
s.t. $A^{(1)}+ A^{(2)}+A^{(3)} = \gamma I$. The space $H$ is
necessarily finite dimensional \cite{ostrovskij1}.  The spectral
decomposition $A^{(i)}=\sum_{k=1}^{m_i} \alpha_k^{(i)} P_k^{(i)}$
with $\alpha_k^{(i)}\not= \alpha_s^{(i)}$ for $k\not=s$ give us the
vectors $\alpha_{i}=(\alpha_k^{(i)})_{k}$ with strictly
 decreasing coefficients and vectors $d^{(i)} = (\rk P_k^{(i)})$. 
To solve the Spectral Problems associated with graph $G$ we need to
describe the set $D = \set{ ( d^{(1)}, d^{(2)}, d^{(3)}, \dim H)}$
corresponding to all irreducible solutions  and for each 
$d\in D$ the set of all possible lists of eigenvalues  $(\alpha^{(1)}, \alpha^{(2)}$, $ \alpha^{(3)}, \gamma)$.

As a result of \cite{KPS} there is an invertible matrix $M_d =
M_d(G)$ s.t. $M_d^{-1}$ is a one-to-one map from  $D$  to a subset
of the  non-degenerate roots of the root system associated with graph $G$
 (non-degenerate roots are difined in Definition~\ref{defnondegen}). We will prove in Section~\ref{s2} that the image of this map is exactly the subset of all non-degenerate roots from the set $\Delta_{sing}\cup\set{\delta_G}$ where $\Delta_{sing}$ is the set
of all singular roots (see Definition~\ref{defnondegen}) and
$\delta_{G}$ is the minimal imaginary root.

The set $\Delta_{sing}$ is a finite union of C-orbits
$\mathcal{C}_1, \ldots, \mathcal{C}_m$ which are  orbits under
the action of Coxeter transformations $\overset{\circ}{c}$ and
$\overset{\bullet}{c}$ ($m$ is the number of vertices in graph $G$).
Each $\mathcal{C}_t$ is of the form $v_s^{(t)} +  \epsilon_t
\mathbb{Z} \delta_G $  where $C_t =\set{v_1^{(t)}, \ldots,
v_{m_t}^{(t)}}$ is a finite set and $\epsilon_t\in\set{1,2}$. These sets are presented in the
appendix. Thus the set $\Delta_{\text{non-degen.}}\cap \Delta_{sing}$ (and hence also $D$) can be explicitly described as a finite union $\Delta_{\text{non-degen.}}\cap\Delta_{sing} = \cup D_t$ of one parameter families of vectors $D_t=\set{d^{(t)}_k, k\ge k_t}$ and  $D = M_d ((\Delta_{\text{non-degen.}}\cap\Delta_{sing})\cup \set{\delta_{G}})$ (matrices $M_d$ for all extended Dynkin graphs are explicitly given in the paper).  For each vector $v= M_d d^{(t)}_k$ we describe the set of all possible eigenvalue vectors $\chi = (\alpha^{(1)}, \alpha^{(2)}, \alpha^{(3)}, \gamma)$ corresponding to $v$ in terms of explicit linear inequalities 
${\rm A}_{t,k}  \chi \ge_t 0$  where matrices ${\rm A}_{t,k}$ are explicitly given 
in the tables in Appendix and for vectors $u$ and $w$ symbol $u\ge_t w$ means that $u_j>w_j$ for $j\not=t$ and $u_t=w_t$. If the vector of multiplicities $v$ correspond to the imaginary root, i.e. $v=M_d \delta_G$ then we give explicit inequalities $H_G$ which define the set of all possible vectors of eigenvalues  $\chi$ using the Horn's inequalities.    
 Thus we completely solve the Spectral problem for all star-shaped extended Dynkin 
graphs.

The final sets of roots $\Delta_f$ for all star-shaped extended Dynkin graph and the orbits of 
the action of the  group generated by the Coxeter transformations $\overset{\bullet}{c}$ and $\overset{\circ}{c}$  on these sets were computed with the help of computer algebra system Wolfram Mathematica.  The Mathematica package can be obtained on  www.math.chalmers.se$\slash\sim$
stanislp$\slash$dynkin$\slash$. The Mathematica package can also be used for computing actions
 of Coxeter transformations on the characters and for computing Horn's inequalities.  
   
The results of this paper were partially announced in~\cite{KPS2}.

\section{
Locally-scalar graph representations. }\label{s1}

The main tools for our classification are Coxeter functors for
locally-scalar graph representations.  First we will recall a
connection between category of $*$-representation of algebra
$\mathcal{A}_{M_1, \ldots, M_n, \gamma}$ associated with a graph
$G$ and the category of locally-scalar representations of the  graph  $G$. For
more details see~\cite{KPS1}.

 For the definitions and notations related to representations of graphs in the category
 of Hilbert spaces we refer the reader to paper~\cite{roiter}.

A graph $G$ consists of a set of vertices $G_v$, a set of edges
$G_e$ and a map $\varepsilon$ from $G_e$ into the set of one- and
two-element subsets of $G_v$ (the edge is mapped into the set of
incident vertices). Henceforth  we consider connected finite
graphs without cycles (trees). Fix a decomposition of $G_v$
 of the form   $G_v= {\overset{\circ}{G}}_v \sqcup
{\overset{\bullet}{G}}_v$ (unique up to permutation) such that for
each  $\alpha\in G_e$ one of the vertices from
$\varepsilon({\alpha})$ belongs to ${\overset{\circ}{G}}_v$ and
the other to ${\overset{\bullet}{G}}_v$.  Vertices in
${\overset{\circ}{G}}_v$ will be called even, and those in the set
${\overset{\bullet}{G}}_v$ odd.
 Let us recall the definition of a  representation  $\Pi$ of a  graph $G$
 in the category of Hilbert spaces  $\mathcal{H}$. Let us associate with each vertex  $g\in G_v$  a Hilbert space
   $\Pi(g)= H_g\in \text{Ob}
\mathcal{H}$, and with each edge  $\gamma\in G_e$ such that
 $\varepsilon(\gamma) =\{g_1,g_2\}$  a pair of mutually adjoint
operators
  $\Pi(\gamma)=\{ \Gamma_{g_1,g_2}$, $\Gamma_{g_2,g_1}\}$, where
  $\Gamma_{g_1,g_2}:H_{g_2}\to H_{g_1}$. We now construct a category
  $\rep (G,\mathcal{H})$.  Its objects are the representations of the
graph  $G$  in  $\mathcal{H}$. A  morphism  $C: \Pi\to
\widetilde{\Pi}$ is a family   $\{C_g\}_{g\in G_v}$ of operators
  $C_g: \Pi(g)\to \widetilde{\Pi}(g)$ such that the following diagrams
  commute  for all edges   $\gamma_{g_2,g_1}\in G_e$:
\begin{equation*}
\begin{CD}
H_{g_1} @>\Gamma_{g_2,g_1}>>H_{g_2} \\  @V C_{g_1} VV  @VV
C_{g_2}V
\\ \widetilde{H}_{g_1}
@>\widetilde{\Gamma}_{g_2,g_1}>>\widetilde{H}_{g_2}
\end{CD}
\end{equation*}

 Let $M_g$ be the set of vertices connected with  $g$ by an edge.
 Let us define the  operators
\[
A_g= \sum_{g'\in M_g} \Gamma_{gg'}\Gamma_{g'g}.
\]

A representation    $\Pi$ in  $\rep (G,\mathcal{H})$  will be
called {\it locally-scalar} if all operators $A_g$ are scalar,
$A_g= \alpha_g I_{H_g}$. The full subcategory of $\rep
(G,\mathcal{H})$,
  the objects  of which are locally-scalar representations, will
be denoted by $\rep G$ and  called the category of locally-scalar
representations of the   graph $G$. 

Let us denote by  $V_{G}$ the real vector space consisting of sequences 
$x=(x_g)$ of real numbers indexed by elements $g\in G_v$. Element $x$ of $V_G$
 will be called a $G$-vector. A vector $x=(x_g)$ is called positive,
$x>0$, if $x\not=0$ and  $x_g\ge 0$ for all $g\in G_v$. Denote
$V_G^+= \{x\in V_G| x>0\}$.  If $\Pi$ is a finite dimensional
representation of the graph $G$ then the  $G$-vector $(d(g))$,
where $d(g)= \dim \Pi(g)$ is called the {\it dimension} of $\Pi$.
If $A_g=f(g) I_{H_g}$ then a $G$-vector $f=(f(g))$ is called a 
{\it character } of the locally-scalar representation   $\Pi$ and
$\Pi$ is called an $f$-representation in this case. The {\it
support} $G_v^{\Pi}$ of
 $\Pi$ is $\{ g\in G_v| \Pi(g) \not= 0 \}$. A representation   $\Pi$
is {\it faithful} if $G_v^{\Pi}= G_v$. A character of the
locally-scalar representation $\Pi$ is uniquely defined on the
support  $G_v^{\Pi}$ and non-uniquely on its complement. In the
general case, denote by $\{f_\Pi \}$ the set of all characters of
$\Pi$.

Below we will describe a connection between locally-scalar graph representation and $*$-representations of algebras $\mathcal{A}_{G, \chi}$ (for details see \cite{KPS1}).

Since we will be concerned with extended Dynkin star-shape graphs we will simplify  
notations and consider only graphs with three rays. This will exclude graph $\widetilde{D}_4$ for
which the formulas are analogous and are left to be recovered by the reader.

 So  we will use notations
$\alpha$, $\beta$, $\delta$ instead of $\alpha^{(1)}$,
$\alpha^{(2)}$, $\alpha^{(3)}$. By $\chi$ we will denote the vector $(\alpha_1,\alpha_2,\ldots, \alpha_k, \beta_1, \beta_2,\ldots,
\beta_l, \delta_1,\delta_2$, $\ldots$, $\delta_m$, $\gamma)$.

Let us consider a tree $G$ with vertices  $\{ g_i,i=0,\ldots,
k+l+m\}$ and edges $\gamma_{g_i g_j}$.

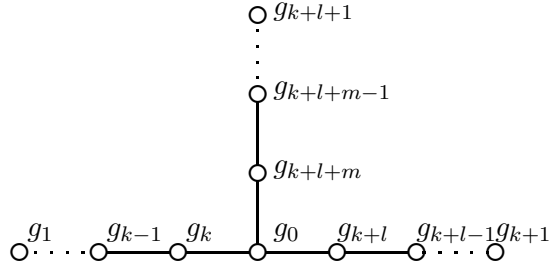
\begin{center}
\setlength{\unitlength}{3pt}
\begin{picture}(70,32)(-1,-1)
\thicklines
 \drawline(29,1)(29,9)
\drawline(29,11)(29,19)
\multiputlist(0,0)(10,0){\circle{2},\circle{2},\circle{2},\circle{2},\circle{2},\circle{2},\circle{2}}
 \multiputlist(30,10)(0,10){\circle{2},\circle{2},\circle{2}}
 \drawline(10,0)(18,0)
\drawline(20,0)(28,0) \drawline(30,0)(38,0)\drawline(40,0)(48,0)
\dottedline{2}(1,0)(7,0) \dottedline{2}(50,0)(58,0)
\dottedline{2}(29,22)(29,28)
 \put(31,10){$g_{k+l+m}$}\put(31,20){$g_{k+l+m-1}$}\put(31,30){$g_{k+l+1}$}
 \put(39,2){$g_{k+l}$}\put(49,2){$g_{k+l-1}$}\put(59,2){$g_{k+1}$}
\put(0,2){$g_{1}$}\put(10,2){$g_{k-1}$}\put(20,2){$g_{k}$}\put(31,2){$g_{0}$}
\end{picture}
\end{center}

Let $\pi$ be a $*$-representation of $\mathcal{A}_{G, \chi}$
 on a Hilbert space  $H_0$. Put
  $P_i=\pi(P_i^{(1)})$, $1\le i \le k$,
 $Q_j=\pi(P_j^{(2)})$, $1\le j \le l$, $S_t=\pi(P_t^{(3)})$, $1\le t \le m$ be the corresponding  spectral projections.

Let $\overline{\rep \mathcal{A}_{G, \chi}}$ denote the full subcategory of ${\rep \mathcal{A}_{G, \chi}}$ consisting of non-degenerate representations. 
In~\cite{KPS1} we constructed a functor which assigns to each representation
$\pi\in \overline{\rep \mathcal{A}_{G, \chi}}$ a locally-scalar representation $\Pi$
 of $G$ with the character  $f$ s.t.  $f(g_0)=\gamma$ and

 \begin{gather}
f(g_k) = \alpha_1,\\\label{beg}
 f(g_{k-1}) =
\alpha_{1}-\alpha_{k},\\
 f(g_{k-2}) = \alpha_{2}-\alpha_{k},\\
f(g_{k-3}) = \alpha_{2}-\alpha_{k-1},\\
f(g_{k-4}) = \alpha_{3}-\alpha_{k-1},\\
\dots\\
f(g_{k+l}) = \beta_1 \\
f(g_{k+l-1}) =
\beta_{1}-\beta_{l},\\
f(g_{k+l+m-2}) =
\beta_{2}-\beta_{m}, \\
f(g_{k+l+m-3}) =
\beta_{2}-\beta_{m-1}, \\
f(g_{k+l+m-4}) =
\beta_{3}-\beta_{m-1},\\
  \dots\\
f(g_{k+l+m}) = \delta_1, \\
f(g_{k+l+m-1})= \delta_{1}-\delta_{m}, \\
f(g_{k+l+m-2}) =
\delta_{2}-\delta_{m}, \\
f(g_{k+l+m-3}) =
\delta_{2}-\delta_{m-1}, \\
f(g_{k+l+m-4}) =
\delta_{3}-\delta_{m-1},\\\label{end}
  \dots
 \end{gather}

If we are given  a locally-scalar representation of the graph
$G$ with a character $f(g_i)=x_i\in\mathbb{R}^*$ then we can construct
a non-degenerate representation of $\mathcal{A}_{G, \chi}$ with  a character
 $(\alpha, \beta, \delta, \gamma)$ s.t.
\begin{gather*}
\alpha_1=x_k,\\
\alpha_k=x_k-x_{k-1},
\alpha_{2}=x_k-x_{k-1}+x_{k-2},\\
\alpha_{k-1}=x_k-x_{k-1}+x_{k-2}-x_{k-3},\\
\alpha_{3}=x_k-x_{k-1}+x_{k-2}-x_{k-3}+x_{k-4},\\
\ldots.
\end{gather*}
Here $x_j=0$ if $j\le 0$. Analogously one can find   $\beta_j$ and
$\delta_t$. We will denote   $\Pi$ by $\Phi(\pi)$.

In fact the above correspondence is given by  a functor   $\Phi:
 \overline{\rep} \mathcal{A}_{G, \chi}\to
 {\rep} (G, f)$, see~\cite{KPS1}.   Moreover, the functor   $\Phi$  is univalent and full.
Let  $\widetilde{\rep}(G,d,f)$ be the full subcategory of
irreducible representations of  ${\rep}(G,d,f)$.
  Let  $\Pi\in \ob \widetilde{\rep}(G,d,f)$, $f(g_i)= x_i\in \mathbb{R}^+,
 d(g_i)=d_i\in \mathbb{N}_0$, where $f$ is  a character of   $\Pi$, $d$ its dimension.
 It easy to verify that the representation $\Pi$ is unitary equivalent to a representation
 from the image of the functor $\Phi$ of the set of non-degenerate representations  if and only if
\begin{gather}\label{nondeg}
 0<x_1<x_2<\ldots < x_k; 0<x_{k+1}<x_{k+2}<\ldots <
x_{k+l};\\
  0<x_{k+l+1}<x_{k+l+2}<\ldots < x_{k+l+m}; \\
0< d_1<d_2<\ldots < d_k< d_0;  0< d_{k+1}<d_{k+2}<\ldots < \label{ineq1}\\
d_{k+l}< d_0;  0< d_{k+l+1}<d_{k+l+2}<\ldots < d_{k+l+m}<
d_0.\label{ineq2}
\end{gather}

  A representation  $\Pi$ of the graph $G$ satisfying conditions
(\ref{nondeg})--(\ref{ineq2})  will be called {\it non-degenerate}.

 Let
\begin{gather*}
\dim H_{p_i}= n_i, 1\le i\le k; \\
\dim H_{q_j}= n_{k+j}, 1\le j\le l; \\
\dim H_{s_t}= n_{k+l+t}, 1\le t\le m; \\
\dim H_{0}= n_{0}.
\end{gather*}
The vector ${n}= (n_0,n_1,\ldots , n_{k+l+m})$ is called the {\it
generalized dimension} of the representation  $\pi$.
Let $\Pi=\Phi(\pi)$ for a non-degenerate
representation $\pi$ of the algebra $\mathcal{A}_{G, \chi}$,
${d}$  = $(d_1,\ldots$, $d_{k+l+m}, d_0)$ be the dimension of  $\Pi$.

The correspondence between generalized dimension of  $\pi$ and that of $\Pi$ is given by the following
equalities:
\begin{eqnarray}
n_0 &=& d_0,\\
n_1+n_2+\ldots + n_k &=& d_k,\\
n_2+\ldots + n_{k-1}+n_k &=& d_{k-1},\\
    n_2+\ldots + n_{k-1} &=& d_{k-2},\\
    n_3+\ldots + n_{k-2}+n_{k-1} &=& d_{k-3},
    \end{eqnarray}
\begin{center}\ldots\end{center}
 Thus
\begin{eqnarray}\label{dim}
n_1 &=&d_k-d_{k-1},\\
n_k &=&d_{k-1}-d_{k-2},\\
n_{2} &=&d_{k-2}-d_{k-3},
\end{eqnarray}
  \begin{center}\ldots\end{center}

  Analogously one can find  $n_{k+1}$, $\ldots$, $n_{k+l}$ from
  $d_{k+1}, \ldots, d_{k+l}$  and  $n_{k+l+1}$, $\ldots$, $n_{k+l+m}$ from  $d_{k+l+1}$, $\ldots$, $d_{k+l+m}$

 Denote by  $\overline{\rep}G$  the full subcategory in   $\rep
 G$ of non-degenerate locally-scalar representations of the graph   $G$.

Now we describe the action of the Coxeter functors on the characters and on the
generalized dimensions of locally-scalar representations. These formulas are necessary for
concrete calculations presented in the last section.

For each vertex $g\in G_v$, denote by  $\sigma_g$ the
linear operator on  $V_G$ given by the formulas:
 \begin{gather*}
 (\sigma_g x)_{g'} = x_{g'}\ \text{if}\ g'\not=g,\\
(\sigma_g x)_{g} = -x_{g} +\sum_{g'\in M_g}x_{g'}.
 \end{gather*}
The mapping   $\sigma_g$  is called the {\it reflection} at the
vertex $g$. The composition of all reflections at odd vertices is
denoted by   $ \overset{\bullet}{c}$ (it does not depend on the
order of the factors), and at all even vertices by $
\overset{\circ}{c}$. A Coxeter transformation is
$c=\overset{\circ}{c}\overset{\bullet}{c}$, $c^{-1}=
\overset{\bullet}{c}\overset{\circ}{c}$.  The transformation
$\overset{\bullet}{c}$ ($\overset{\circ}{c}$)  is called an  odd
(even) Coxeter map.
  Let us adopt the  following notations for compositions of the  Coxeter maps:
   $\overset{\bullet}{c}_k= \ldots
\overset{\bullet}{c}\overset{\circ}{c}\overset{\bullet}{c}$ ($k$
factors), $\overset{\circ}{c}_k= \ldots
\overset{\circ}{c}\overset{\bullet}{c}\overset{\circ}{c}$ ($k$
factors), $k\in\mathbb{N}$.

Any real function  $f$ on $G_v$ can be identified  with a
$G$-vector $f=(f(g))_{g\in G_v}$. If $d(g)$  is the dimension of a
locally-scalar graph representation  $\Pi$, then
\begin{gather}
\overset{\circ}{c}(d)(g)=\begin{cases} -d(g) +  \sum_{g'\in
M_g}d(g'), &\text{if}\  g\in \overset{\circ}{G}_v,\\
d(g), & \text{if}\ g\in \overset{\bullet}{G}_v,
\end{cases}\\
\overset{\bullet}{c}(d)(g)=\begin{cases} -d(g) +  \sum_{g'\in
M_g}d(g'), &\text{if}\  g\in \overset{\bullet}{G}_v,\\
d(g), & \text{if}\ g\in \overset{\circ}{G}_v.
\end{cases}
\end{gather}

 For $ d\in Z_G^+$ and  $f\in V_G^+$, consider the full subcategory
  $\rep (G,d,f)$ in $\rep
 G$ (here $Z_G^+$ is the set of positive integer  $G$-vectors),
with the set of objects  $Ob \rep (G, d, f)= \{ \Pi| \dim \Pi (g)
=d(g), f\in\{f_\Pi\}\}$. All representations
    $\Pi$ from $\rep (G,d,f)$ have the same support  $X=X_d=G_v^\Pi=\{g\in G_v|
  d(g)\not=0\}$.   We will consider these categories only if   $(d,f)\in
  S=\{ (d,f)\in Z_G^+\times V_G^+ | d(g) + f(g) >0, g\in G_v \}$.
   Let  $\overset{\circ}{X}= X\cap \overset{\circ}{G}_v$,
    $\overset{\bullet}{X}= X\cap \overset{\bullet}{G}_v$.  $\rep_{\circ}(G, d, f) \subset \rep(G, d,
    f)$  ( $\rep_{\bullet}(G, d, f) \subset \rep(G, d,
    f)$) is the full subcategory with objects  $(\Pi, f)$ where
      $f(g)>0$  if  $g\in \overset{\circ}{X}$ ($f(g)>0$  if  $g\in
    \overset{\bullet}{X}$). Let $S_0=\{(d,f)\in S|f(g)>0 \text{ if }
g\in \overset{\circ}{X}_d\}$, $S_\bullet=\{(d,f)\in S|f(g)>0
\text{ if } g\in \overset{\bullet}{X}_d\}$

Put
\begin{gather}
\overset{\bullet}{c}_d(f)(g)=\overset{\circ}{f}_d(g)=
\begin{cases}
 \overset{\bullet}{c}(f)(g), &\text{if}\  g\in \overset{\bullet}{X}_d,\\
f(g), & \text{if}\ g\not\in \overset{\bullet}{X}_d,
\end{cases}\\
\overset{\circ}{c}_d(f)(g)=\overset{\bullet}{f}_d(g)=
\begin{cases}
 \overset{\circ}{c}(f)(g), &\text{if}\  g\in \overset{\circ}{X}_d,\\
f(g), & \text{if}\ g\not\in \overset{\circ}{X}_d.
\end{cases}.
\end{gather}

 Let us denote  $$\overset{\bullet}{c}_d^{(k)}(f)=
 \ldots \overset{\bullet}{c}_{\overset{\circ}{c}_2(d)}
 \overset{\circ}{c}_{\overset{\circ}{c}(d)}
 \overset{\bullet}{c}_d(f) \quad (k \text{ factors}) $$

 $$\overset{\circ}{c}_d^{(k)}(f)=
 \ldots \overset{\circ}{c}_{\overset{\bullet}{c}_2(d)}
  \overset{\bullet}{c}_{\overset{\bullet}{c}(d)}
  \overset{\circ}{c}_d(f) \quad (k \text{ factors}) $$
The even and odd Coxeter reflection functors are defined
in~\cite{roiter}. They act between categories as follows: 
  $$\overset{\circ}{F}: \rep_{\circ} (G, d, f) \to
  \rep_{\circ} (G, \overset{\circ}{c}(d), \overset{\circ}{f}_d)\quad \text{ if }  (d,f)\in S_\circ $$

$$\overset{\bullet}{F}: \rep_{\bullet} (G, d, f) \to
  \rep_{\bullet} (G, \overset{\bullet}{c}(d), \overset{\bullet}{f}_d)\quad \text{ if } (d,f)\in S_\bullet$$
   These functors  are equivalences of the  categories.
Let us denote
   $\overset{\circ}{F}_k(\Pi)=
 \ldots \overset{\circ}{F}\overset{\bullet}{F}\overset{\circ}{F}(\Pi)$ ($k$  factors),
 $\overset{\bullet}{F}_k(\Pi)=
 \ldots \overset{\bullet}{F}\overset{\circ}{F}\overset{\bullet}{F}(\Pi)$ ($k$
 factors),  if the compositions make sense. 
 
 Using these functors, an analog of Gabriel's theorem
for graphs and their locally-scalar representations has been
proven in~\cite{roiter}. In particular, it has been proven that
any locally-scalar graph representation of a Dynkin graph decomposes into a direct
sum (finite or infinite) of finite-dimensional indecomposable
representations, and all indecomposable representations can be
obtained by odd and even Coxeter reflection functors starting from
the simplest representations  $\Pi_g$ of the graph  $G$. The simplest representation is the one which has coordinate vector as its generalized dimension, i.e. it corresponds to a vertix $g$:  $\Pi_g(g)=\mathbb{C}, \Pi_g(g')= 0$  if $g\not=g';\ g,g'\in G_v
 $.

In sequel we will refer to Coxeter functors between the categories of $*$-representations of algebras $\mathcal{A}_{G, \chi}$. We define these functors simply as  $S=
\Phi^{-1}\overset{\circ}{F} \Phi $ and $T=
\Phi^{-1}\overset{\bullet}{F} \Phi$.

\section{Root systems for extended Dynkin graphs.}\label{s2}

 Let us recall a few facts about root systems associated with extended Dynkin graphs (see  \cite{cw}, \cite{redchuk}).
 Let $G$ be a simple connected graph. Then its {\it Tits form} is the following quadratic form
 $$q(\alpha)= \sum_{i\in G_v} \alpha_i^2 -\frac{1}{2} 
\sum_{} \alpha_i\alpha_j \quad (\alpha \in V_G)
 $$
 where the second summation is over pairs $(i,j)$ such that there exists edge 
 $\beta\in G_e$ with $\set{i,j}=\epsilon(\beta)$, i.e. each edge is counted twice.  
The corresponding {\it symmetric bilinear form} is  $(\alpha, \beta) = q(\alpha +
\beta ) - q(\alpha) -q(\beta)$. Vector $\alpha \in V_G$ is called
 {\it sincere} if each component is non-zero.

 It is well known that for Dynkin graphs (and only for them)
 bilinear form $(\cdot, \cdot)$ is positive definite. The form is
 positive semi-definite for extended Dynkin graphs. And in the
 letter case $\rad q= \{v | q(v)=0 \}$ is equal to $\mathbb{Z}
 \delta$ where  $\delta$  is a minimal imaginary root. For other
 graphs (which are neither Dynkin nor extended Dynkin) there are
 vectors $\alpha \ge 0$ such that $q(\alpha)<0$ and $(\alpha, \epsilon_j)\le
 0$ for all $j$.

 For an extended Dynkin graph $G$ a vertex $j$ is called
 {\it extending } if $\delta_j =1$. The graph obtained by deleting
 extending vertex is the {\it corresponding} Dynkin graph. The
 {\it set of roots} is $\Delta = \{ \alpha \in V_G | \alpha_i \in \mathbb{Z} \text{ for all } i\in G_v, \alpha\not=0, q(\alpha)\le 1
 \}$. A root $\alpha$ is {\it real} if $q(\alpha) = 1$ and
 {\it imaginary} if $q(\alpha) =0$. Every root is either positive
 or negative, i.e. all coordinates are simultaneously non-negative or non-positive.

\begin{definition}\label{defnondegen}
\
\begin{enumerate}
\item Root $d$ satisfying \eqref{ineq1}-\eqref{ineq2} will be called
non-degenerate.

\item 
Root $d$ is called regular if $c^t(d)\in V_G^+$ for every $t\in\mathbb{Z}$ and singular otherwise. Here $c=\overset{\circ}{c}\overset{\bullet}{c}$ is a Coxeter transformation.   

\end{enumerate}
\end{definition}
The definition of regular and singular roots are due to Gelfand and Ponomarev~\cite{GP}.  

By proposition \cite[Proposition~1.6]{GP} for any extended Dynkin graph $G$  root $v$ is regular if and only if  and only if  $L_G(v)=0$ where $$L_G(v) =\sum_{g\in \overset{\bullet}{G}} u_g v_g- \sum_{h\in \overset{\circ}{G}}u_h v_h$$ and $\delta_G = (u_g)_g$ is the minimal imaginary root. The root $w$ is singular if and only if it belongs to the orbit of some coordinate vector $\varepsilon_j$ under the action of the group generated by Coxeter transformations $\overset{\bullet}{c}$ and $\overset{\circ}{c}$ see~\cite[Proposition~1.4]{GP}.

Since by \cite[Lemma~1.2]{GP} for every $x\in V_g$,  $L_G(\overset{\bullet}{c}(v)) = L_G(\overset{\circ}{c}(v)) = - L_G(v)$  the above mentioned characterization of the regular root implies that the set of all regular roots  $\Delta_{reg}$ is invariant with respect to Coxeter maps $\overset{\circ}{c}$ and $\overset{\bullet}{c}$.

 For our classification purposes we will need the
 following fact~\cite[p.18]{cw}:
\begin{statement} For an extended Dynkin graph the
 set $\Delta\cup \{0\}\slash \mathbb{Z} \delta$ is finite. Moreover, if $e$ is an
 extending vertex then the set $\Delta_f =\{\alpha \in \Delta\cup\{0\} |
 \alpha_e=0 \}$ is a complete set of representatives of the cosets
 from $\Delta\cup \{0\} \slash \mathbb{Z} \delta$.
 \end{statement}

If $\alpha$ is a root then $\alpha+\delta$ is again a root.
\begin{definition}
We call the coset $\alpha +\delta\mathbb{Z}$ the
$\delta$-series and the coset $\alpha +2\delta\mathbb{Z}$ the $2\delta$-series.
 If $\alpha$ is a root then its images under the action  of the group generated by
 $\overset{\circ}{c}$ and $\overset{\bullet}{c}$ will be called a Coxeter series or
 $C$-series for short.
 \end{definition}
 
 As a result of calculations presented in the tables at the end of the paper we have that   $C$-series of any singular root decomposes into a finite number of $\delta$-series or $2\delta$-series of roots.

Note that  to find  formulas of the locally-scalar representations
of a given extended Dynkin graph we need to consider two
principally different cases: the case when the vector of
generalized dimension is a real root and the case when it is an
imaginary root. In the letter case the vector of parameters $\chi$
must satisfy  a certain
linear relation which is obtained by taking traces from the both sides of the equation $A_1+\ldots +A_n =\gamma I$.  Hence $\chi$ must belong to a
certain hyperplane $h_G$ which depends only on the graph $G$. A
simple calculation using \eqref{beg}-\eqref{end} yields that for extended Dynkin graphs
$\widetilde{D}_4$, $\widetilde{E}_6$, $\widetilde{E}_7$,
$\widetilde{E}_8$ these hyperplanes are the following:

\begin{description}
\item[$\widetilde{D}_4$:] $\alpha_1+\beta_1+\delta_1 +
\eta_1 = 2 \gamma$
\item[$\widetilde{E}_6$:] $\alpha_1+\alpha_2+\beta_1+\beta_2
+\delta_1+\delta_2 = 3 \gamma$
\item[$\widetilde{E}_7$:] $\alpha_1+\alpha_2+\alpha_3+
\beta_1+\beta_2 +\beta_3+ 2\delta_1 = 4 \gamma$
\item[$\widetilde{E}_8$:] $\alpha_1+\alpha_2+\alpha_3+\alpha_4+ 2(\beta_1+\beta_2)  +3\delta_1 = 6 \gamma$
\end{description}

It is know (see~\cite{melit}) that in case $\chi \in h_G$ the
dimension of any irreducible representation is bounded (by 2 for
$\widetilde{D}_4$, by 3 for $\widetilde{E}_6$, by 4 for
$\widetilde{E}_7$ and by 6 for $\widetilde{E}_8$). Thus in case $\chi \in h_G$ we can describe the set of admissible parameters
$\chi$ using  Horn's inequalities. In case $\chi\not\in h_G$
 the dimension of any  irreducible locally-scalar representation
 is a real root. We will show that only singular roots can occur. First we need some
 definitions from~\cite{ostrovskij}. To avoid conflict of terminology we will call
 vector $
\chi=(\alpha_1^{(1)},\dots,\alpha_{m_1}^{(1)};\dots;
\alpha_1^{(n)},\dots,\alpha_{m_n}^{(n)})$ a {\it reduced character} of the graph. Clearly the character is obtained from
reduced character by appending $\gamma$.

We will recall below the formulas for the action of the Coxeter functors on the characters.
Recall that
$S$ and $T$ are the Coxeter functors defined in the previous section.  Then
 $S\colon \rep \mathcal
A_{G,\chi,\lambda} \to \rep \mathcal A_{G,\chi',\lambda'}$
and $T\colon \rep \mathcal
A_{G,\chi,\lambda} \to \rep \mathcal
A_{G,\chi'',\lambda}$, where
\begin{align*}
\chi'&=(\alpha_{1}^{(1)}-\alpha_{0}^{(1)},
 \alpha_{1}^{(1)}-\alpha_{m_1}^{(1)}, \dots,
\alpha_{1}^{(1)}-\alpha_{2}^{(1)}; \dots; \alpha_{1}^{(n)}-
\alpha_{0}^{(n)}, \alpha_{1}^{(n)}-
\alpha_{m_n}^{(n)},\\
 &\dots, \alpha_{1}^{(n)}-\alpha_{2}^{(n)}),
\\
\lambda' &= \alpha_{1}^{(1)} + \dots +\alpha_{1}^{(n)} -\lambda,
\\
\chi''&=  (\lambda-\alpha_{m_1}^{(1)},\dots, \lambda -\alpha_1^{(1)};
\dots ;\lambda-\alpha_{m_n}^{(n)},\dots, \lambda -\alpha_1^{(n)}).
\end{align*}
The action of these functors on $*$-representations gives rise to the
action on pairs, $S\colon (\chi,\lambda) \mapsto (\chi',\lambda')$,
$T\colon (\chi,\lambda )\mapsto (\chi'',\lambda)$.

  Let $\chi$ be a reduced  character on the graph $G$,
\begin{gather}\label{chi}
\chi = (\alpha_1^{(1)}, \dots, \alpha_{m_1}^{(1)}; \dots;
\alpha_1^{(n)},
\dots,\alpha_{m_1}^{(n)}),
\\
0<\alpha_{m_l}^{(l)}<\dots
<\alpha_{1}^{(l)},\quad
l=1, \dots, n. \notag
\end{gather}

The following notion of invariant functional was introduced in~\cite{ostrovskij}.
Let $\omega(\cdot)$ be a linear functional, which takes non-negative
values on reduced characters.

\begin{definition}
We say that $\omega(\cdot)$ is invariant with respect to the functor
$TS$, if
\[
TS(\chi,\omega(\chi)) = (\tilde\chi,\omega(\tilde\chi))
\]
for any reduced character $\chi$ on $\Gamma$.
\end{definition}

If  $G$ is an extended Dynkin Graph then there exists unique
   invariant functional:

--- for $\tilde D_4$, $\omega(\alpha;\beta;\gamma;\delta) =
    \frac12(\alpha+\beta +\gamma +\delta)$;

--- for $\tilde E_6$,
    $\omega(\alpha_1,\alpha_2;\beta_1,\beta_2;\gamma_1,\gamma_2) =
    \frac13(\alpha_1 +\alpha_2 +\beta_1 +\beta_2 +\gamma_1 +\gamma_2)$;

--- For $\tilde E_7$,
    $\omega(\alpha_1,\alpha_2,\alpha_3;\beta_1,\beta_2
    ,\beta_3;\gamma) = \frac14(\alpha_1 +\alpha_2 +\alpha_3 +\beta_1 +\beta_2
    +\beta_3 +2\gamma)$;

--- for $\tilde E_8$, $\omega(\alpha_1,\alpha_2,
    \alpha_3,\alpha_4,\alpha_5;\beta_1,\beta_2;\gamma) =  \frac16(\alpha_1
    +\alpha_2 +\alpha_3+\alpha_4 + \alpha_5 + 2 \beta_1 +2\beta_2
    +3\gamma)$.

Note that the equation $\lambda = \omega(\chi)$ defines the hyperplane $h_{G}$.

\begin{theorem}\label{ostrovth}
Let $\pi$ be an irreducible non-degenerate $*$-representation of the algebra
$\mathcal{A}_{G,\chi, \lambda}$ associated with an extended Dynkin graph $G$
and $\widehat{\pi}$ be the corresponding locally-scalar representation of the graph $G$.
Then either generalized dimension $d$ of $\widehat{\pi}$ is a
singular root or  $(\chi, \lambda)\in h_G$.
\end{theorem}
\begin{proof}
Let $\pi$ be irreducible non-degenerate representation of
$\mathcal{A}_{G,\chi, \lambda}$ and $\hat{\pi}= \Phi(\pi)$ be the
corresponding locally-scalar representation
of $G$. Then the vector of generalized dimension $d$ of $\hat{\pi}$
is non-degenerate. Suppose that $d$ is regular. 
Assume  that $\lambda < \omega(\chi)$. By the proof of Theorem 3
in~\cite{ostrovskij} there exists positive integer $n$ s.t. representation  $\pi_n = (ST)^n(\pi)$ corresponds to  a proper subgraph and thus $\Phi(\pi_n)$ is not non-degenerate. Hence its generalized dimension is a singular root. Since the generalized dimension of $\Phi(\pi_n)$ and $d$ clearly belong to the same $C$-series  by invariance of the set of singular (regular) roots  with respect to Coxeter maps  we have that $\Phi(\pi_n)$ must be
singular. This contradiction proves that $\lambda \ge \omega(\chi)$.
If $\lambda >\omega(\chi)$ we can apply  functor $S$  to $\pi$ and get representation
$\pi'$ corresponding to a pair $(\chi',\lambda')$ with
   $\lambda'<\omega(\chi')$. The vector of generalized dimension of
   $\Phi(\pi')$ is non-degenerate and regular and by the first part of the proof we get a contradiction. Thus
   $\lambda = \omega(\chi)$ and hence $(\chi, \lambda)\in h_G$.
\end{proof}

\begin{lemma}\label{regl}

\begin{enumerate}\label{lemreg}
The following statements hold
\item
The minimal dimension of the non-degenerate locally-scalar  representation of
an extended Dynkin graph with vector of generalized dimension being regular root is not less than 

--For $\widetilde{D}_4$, $3$;

--For $\widetilde{E}_6$, $4$;

--For $\widetilde{E}_7$, $5$;

--For $\widetilde{E}_8$, $7$.
\item
There is no irreducible non-degenerate representation of $\mathcal{A}_{G, \chi}$ for
$\chi\in h_G$ in dimension  $M_d w$ where
$w$ is a real regular root for any extended Dynkin graph $G$.
\item
There are no irreducible non-degenerate  $*$-representations of
the algebra  $\mathcal{A}_{G, \chi}$ of generalized dimension $M_d w$
where $w$ is  a real regular root.
\end{enumerate}
\end{lemma}
\begin{proof}

For any root $v$ denote by $ng(v)$ the non-degenerate root of the form  $v+ k \delta$ with minimal possible
integer $k\ge 0$ (it is clear that such $k$ exists). Note that if $v$ is
non-degenerate then clearly the same is true for $v+k \delta$ for any $k\ge 0$.
It is also clear that any non-degenerate regular $v$ belongs to $ng(\Delta_{reg}
\cap \Delta_{f}) + \mathbb{Z}_+ \delta$. 

Using the tables 7-9 of the regular roots belonging to  $\Delta_{reg}\cap \Delta_f$ and computing $ng(v)$ for all $v\in \Delta_{reg}\cap \Delta_f$. We obtain a final sets of roots for each extended Dynkin graph $G$ and   the minimum of the last coordinates of these vectors will be $3,4,5,7$ for $G= \widetilde{D}_4, \widetilde{E}_6, \widetilde{E}_7, \widetilde{E}_8$ correspondingly which proves the first statement.

To prove the second one recall that
$\mathcal{A}_{G, \chi}$ is $\mathbb{F}_{2n}$-algebra for any $\chi\in h_{G}$
with $n= 4$ for $\widetilde{D}_4$, $n=9$ for $\widetilde{E}_6$, $n = 16$ for
$\widetilde{E}_7$ and  $n=36$ for $\widetilde{E}_8$
(see~\cite[Theorem~3]{melit1}). Hence the maximal dimension of the irreducible
representation of $\mathcal{A}_{G, \chi}$ is 2 for
$\widetilde{D}_4$,  3 for $\widetilde{E}_6$,  4 for
$\widetilde{E}_7$ and  6 for $\widetilde{E}_8$. Since the  dimension of the algebra
representation is the same as the dimension  of the
corresponding locally-scalar representation of the graph
the second claim follows from the first one.

For $\chi\not\in h_G$ the third claim  follows
from Theorem~\ref{ostrovth}. For $\chi\in h_G$  the claim follows from
the second claim.
\end{proof}

For a vector $v=(v_0, \ldots, v_n)$ and $0\le s\le n$ we will write
$v\ge_s 0$ if $v_j>0$ for all $j\not=s$ and $v_s=0$.

The equivalence functor $\Phi$  assigns to every representation
$\pi\in \mathcal{A}_{G, \chi}$ of generalized dimension $(l_1,\ldots, l_n)$
a unique  locally-scalar representation of graph $G$ with a
character $(x_1, \ldots, x_n, x_0)$ and dimension
$(v_1, \ldots, v_n, v_0)$.
Let $M_f$ denote the transition matrix which transform the
vector $\chi$ to $(x_1, \ldots, x_n, x_0)$, i.e. $M_f \chi^t = (x_1, \ldots, x_n, x_0)^t$
(where $v^t$ denote the
transposed vector $v$).
Let  $M_d$ be the transition matrix which transforms generalized
dimension $(v_1, \ldots, v_n, v_0)$  of a graph  representation
 to generalized dimension $(l_1,\ldots, l_n)$ of the
 corresponding algebra representation, i.e. $M_d (v_1, \ldots, v_n, v_0)^t = (l_1,\ldots, l_n)^t $.
 Further we will omit $t$ superscript and write $M_f v$ instead of $M_f v^t$.

In the following theorem we present the general form of a solution of Spectral problem for
 extended Dynkin graphs. In a subsequent sections we present an explicit solution separately
 for each graph.

\begin{theorem}\label{genform} Let $G$ be an extended Dynkin graph and $\pi$ be a non-degenerate
irreducible $*$-representation in a generalized dimension $v$ of
$\mathcal{A}_{G, \chi}$ for some character $\chi$. Then one  of two possibilities
holds
\begin{itemize}
\item $\chi \in h_G$ and $d= M_d\delta$ where $\delta$ is the minimal imaginary
root of the root system associated with
$G$.
\item There $k$ and $t$ such that

\begin{gather}
\overset{\bullet}{c}^{(k)}_d M_f \chi \ge_t 0,\\
\overset{\bullet}{c}_{k}(M_d^{-1}v)=e_t,
\end{gather}
or
\begin{gather}
\overset{\circ}{c}^{(k)}_d M_f \chi  \ge_t 0,\\
\overset{\circ}{c}_{k}(M_d^{-1}v)=e_t.
\end{gather}
(depending on the parity of $k+t$). Moreover, systems of
inequalities \eqref{cond1}, \eqref{cond2}
 are necessary and sufficient conditions for the
existence of representation of $\mathcal{A}_{G, \chi}$ in dimension $v$.
\end{itemize}
\end{theorem}
\begin{proof}
By Theorem~\ref{ostrovth} and Lemma~\ref{lemreg} vector $v$ is of the form $M_d w$ where $w$ is either minimal imaginary
root $\delta$ or $w$ is a singular root.
If $\Pi\in \rep (G,w,\xi)$ is a locally-scalar representation of a
graph $G$ with $w$ being a singular root
then $e_t= \overset{\bullet}{c}_k(w)$ or $e_t= \overset{\circ}{c}_k(w)$
for some positive integer $k$ and a coordinate vector $e_t$.
Thus there is a locally-scalar representation $\Pi' \in \rep (G,e_t,\xi')$
such that applying corresponding sequence of
Coxeter functors $\ldots\overset{\circ}{F}\overset{\bullet}{F}$
or $\ldots\overset{\bullet}{F}\overset{\circ}{F}$ to $\Pi'$ we obtain
$\Pi$ and  hence $w$ belong to $C$-orbit of $e_t$ and
$\overset{\bullet}{c}^{(k)}_d(\xi) = \xi'$ or
$\overset{\circ}{c}^{(k)}_d(\xi) = \xi'$.
Thus the necessary and sufficient conditions on $\xi$ for
the representation $\Pi$ to exist can be written as
inequalities \begin{gather}
\overset{\bullet}{c}^{(k)}_d(\xi) \ge_t 0,\label{cond1}
\end{gather}
or
\begin{gather}
\overset{\circ}{c}^{(k)}_d(\xi) \ge_t 0, \label{cond2}
\end{gather}

Let $C_t$ denote the $C$-series of $e_t$. It can be checked by
direct computations that for each extended Dynkin graph every $C$-series
  is a union of finite number of $\delta$-series or $2\delta$-series of roots,
   i.e. $C_t=(v_0+\epsilon \delta \mathbb{Z})\cup \ldots \cup (v_m+\epsilon
   \delta \mathbb{Z})$ where $\epsilon\in\{1, 2\}$
   and $\overset{\circ}{c}(v_{2r})= v_{2r+1}$,
   $\overset{\bullet}{c}(v_{2r-1})= v_{2r}$
   (or $\overset{\bullet}{c}(v_{2r})= v_{2r+1}$,
   $\overset{\circ}{c}(v_{2r-1})= v_{2r}$).
   We have presented this finite sequences
   $(v_0, \ldots, v_m)$ in the tables at the end of the paper.
   Elements of $C$-series  can be written then as
   $w_k = v_{k\mod (m+1)} + \epsilon [\frac{k}{m+1}]\delta$
   where $[x]$ denote the integer part of $x$.

 Then the generalized
 dimension $l$ of the  irreducible  representation $\pi\in \rep \mathcal{A}_{G, \chi}$ is of the form $M_d w$ where $w$ in a non-degenerate root of root system associated with graph $G$ and conditions (\ref{cond1}), (\ref{cond2}) give the following necessary and sufficient conditions on $\chi$ of existence of representation in dimension $l$:
$$\overset{\bullet}{c}^{(k)}_d M_f \chi \ge_t 0$$ or
$$\overset{\circ}{c}^{(k)}_d M_f \chi  \ge_t 0.$$

\end{proof}

Let $t$ be the minimal number  such that $w_t$ is non-degenerate then
$w_l$ is also non-degenerate for all $l>t$.
We will denote by $D_{j,t}$ the matrix which transform
the character of a locally-scalar graph representation with
dimension $w_t$ to the one with dimension $e_j$, i.e.
 $D_{j, k} (x_1,x_2,\ldots,x_n, x_0)^t =
 (x_1',x_2',\ldots,x_7', x_0')^t$ where
 $(x_1',x_2',\ldots,x_7', x_0')$ obtained
 from  $(x_1,x_2$, $\ldots$, $x_7, x_0)$
 by applying the corresponding  sequence of
 Coxeter maps that transform $v_k$ to $e_j$.

The explicit results of the computations of $C$-series,
matrices $M_f$, $M_d$ etc. are gathered in the tables in Appendix.
In the following sections we present explicit answer to the Spectral problem
for each extended Dynkin graphs separately.

\section{ Representations of  $\mathcal{A}_{\widetilde{D}_4,\chi}$. }

The parameters $\chi=(\alpha,   \beta, \xi,  \delta,   \gamma)$ of the algebra $\mathcal{A}_{\widetilde{D}_4,\chi}$ and the vector of generalized dimension $n=(n_1, \ldots, n_4, n_0)$ will be plotted on the associated graph according to the following picture:
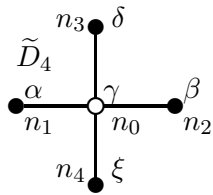
\begin{center}
\setlength{\unitlength}{3pt}
\begin{picture}(30,30)(-1,-1)
\thicklines
\put(-1,15){$\widetilde{D}_4$}
 \drawline(9,11)(9,19)
 \drawline(9,9)(9,1)\put(9,0){\circle*{2}}
 \multiputlist(0,10)(10,0){\circle*{2},\circle{2},\circle*{2}}
 \multiputlist(10,20)(0,10){\circle*{2}}
 \drawline(0,10)(8,10)
\drawline(10,10)(18,10)
 \put(0,11){$\alpha$}
\put(10,11){$\gamma$}
 \put(20,11){$\beta$}
\put(11,20){$\delta$}
\put(11,1){$\xi$}
\put(0,7){$n_1$}
\put(11,7){$n_0$}
 \put(20,7){$n_2$}
\put(4,20){$n_3$}
\put(4,1){$n_4$}
\end{picture}
\end{center}
The category of non-degenerate $*$-representations of $\mathcal{A}_{\widetilde{D}_4,\chi}$ is equivalent to the category of non-degenerate locally-scalar representations of the graph $\widetilde{D}_4$ with the character and generalized dimension given on the following picture:

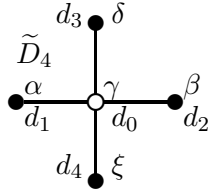
\begin{center}
 \setlength{\unitlength}{3pt}
\begin{picture}(30,30)(-1,-1)
\thicklines
\put(-1,15){$\widetilde{D}_4$}
 \drawline(9,11)(9,19)
 \drawline(9,9)(9,1)\put(9,0){\circle*{2}}
 \multiputlist(0,10)(10,0){\circle*{2},\circle{2},\circle*{2}}
 \multiputlist(10,20)(0,10){\circle*{2}}
 \drawline(0,10)(8,10)
\drawline(10,10)(18,10)
 \put(0,11){$\alpha$}
\put(10,11){$\gamma$}
 \put(20,11){$\beta$}
\put(11,20){$\delta$}
\put(11,1){$\xi$}
\put(0,7){$d_1$}
\put(11,7){$d_0$}
 \put(20,7){$d_2$}
\put(4,20){$d_3$}
\put(4,1){$d_4$}
\end{picture}
\end{center}

Obviously the transition matrix  $M_f$ such that $M_f(\chi)=(x_1,x_2,\ldots, x_4,x_0)$ and  the transition matrix $M_d$ such that $M_d(d_1,\ldots, d_4, d_0)^T = (n_1, \ldots, n_4, n_0)^T$ are identity matrices.

For  $\Pi\in {\rep}(G,d,f)$ with sincere $d$ the Coxeter map $\overset{\bullet}{C}\overset{\circ}{C}$ transform character $x= (x_1,\ldots,x_4, x_0)$ by multiplying from the left the vector-column $x$ on the matrix $M_c= {\tiny \left(%
\begin{array}{ccccc}
  -1 & 0 & 0 & 0 & 1 \\
 0 & -1 & 0 & 0 & 1 \\
 0 & 0 & -1 & 0 & 1 \\
 0 & 0 & 0 & -1 & 1 \\
 -1 & -1 & -1 & -1 & 3 \\
\end{array}%
\right).}$

It is easy to check that  the Jordan form of $M_c$ is
\[J = (-1)\oplus (-1)\oplus (-1) \oplus J_2(1) ,\] where \[J_2(1)= \left(%
\begin{array}{cc}
  1 & 1 \\
  0 & 1 \\
\end{array}%
\right).\]

Specializing Theorem~\ref{genform} to the case  $G= \widetilde{D}_4$ we get  Theorem~\ref{ineqeD4} below. In view of the symmetry of the graph $\widetilde{D}_4$ the $C$-series and  $\delta$-series of the coordinate vectors corresponding to vertices $2, 3, 4$ differ from the $C$-series and $\delta$-series of the coordinate vector corresponding to vertex $1$ by corresponding transpositions $\tau_1=(1,2), \tau_2=(1,3), \tau_3=(1,4)$ of coordinates. 
We will write $\tau v$ to denote vector obtained by permuting coordinates of vector $v$ according to permutations $\tau$. Only  the $C$-series of the roots $(1,0,0,0,0)$ and $(0,0,0,0,1)$ are listed in the appendix. The rest of  $C$-series can be obtained by applying permutations from the set $\set{(1,2),  (1,3), (1,4)}$.    The notation $u \ge_t v$ used in theorem below has been introduced  before Theorem~\ref{genform}. 

Vectors $d^{(t)}_1, d^{(t)}_2, \ldots, $ represent  roots from the $C$-series of coordinate vector $\varepsilon_t$ corresponding to the vertex $t$ of the graph listed  in the same  order they appear in the orbit of the Coxeter transformation $c$. Number $k_t$ denotes the minimal number of Coxeter transformations $\overset{\bullet}{c}$ and $\overset{\circ}{c}$ necessary to apply (in alternating order) to get from coordinate vector $\varepsilon_t$ to a non-degenerate root.   

\begin{theorem}\label{ineqeD4}
Let $d^{(t)}_k = v^{(t)}_{k \mod m_t}+
[\frac{k}{m_t}]\delta$ for $k\ge k_t$ and $t\in\{0,1\}$. Here  $v^{(t)}_s$ is the
$s$-th vector in the set $C_t$  from Table~\ref{tabD4}, $m_t = |C_t|$ and $k_1 = 5$, $k_0=2$.

If the vector $\chi\not\in h_{\widetilde{D}_4}$, i.e. $\chi$ does not satisfy $\alpha_1+\beta_1+\delta_1 +
\eta_1 = 2 \gamma$  then
 the algebra  $\mathcal{A}_{\widetilde{D}_4,\chi}$ has an irreducible non-degenerate representation   in a generalized dimensions $v$ if and only if two conditions holds 
 \begin{enumerate} \item for some $t\in\{0,1\}$,  $k\ge
k_t$ and a permutation $\tau\in \{e,   (1,2),  (1,3), (1,4)\}$ ($\tau$ is identity $e$ if $t=0$), $$v=  \tau M_d d^{(t)}_k.$$ 

\item $${\rm A}_{t,k}\tau \chi \ge_t 0$$ where the matrix ${\rm A}_{t,k}$ is
taken from Table~\ref{tabD4}. 
\end{enumerate}
The permutation $\tau$ and numbers $t$ and  $k$ are determined uniquely.    The irreducible representation of the algebra $ \mathcal{A}_{\widetilde{D}_4,\chi}$ corresponding to the above generalized dimension and $\chi$ is unique.

If $\chi$  satisfies the equation $\alpha_1+\beta_1+\delta_1 +
\eta_1 = 2 \gamma$ then irreducible non-degenerate
representations of $\mathcal{A}_{\widetilde{D}_4,\chi}$ may exist only
in the generalized dimension $\delta_{alg}(\widetilde{D}_4)$.
They exist if and only if $\chi$ satisfies inequalities $H_{\widetilde{D}_4}$ from Table~\ref{tabD4}.    \end{theorem}

\section{ Representations of  $\mathcal{A}_{\widetilde{E}_6,\chi}$. }

The parameters $\chi=(\alpha_1, \alpha_2,  \beta_1, \beta_2,  \delta_1, \delta_2,  \gamma)$ of the algebra $\mathcal{A}_{\widetilde{E}_6,\chi}$ and the vector of generalized dimension $d=(d_1, \ldots, d_6, d_0)$ will be plotted on the associated graph according to the following picture:
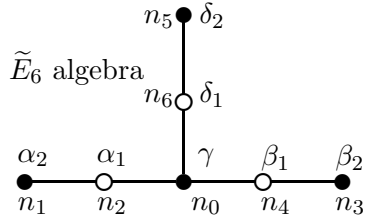
\begin{center}
\setlength{\unitlength}{3pt}
\begin{picture}(50,30)(-1,-7)
\thicklines
\put(-3,13){$\widetilde{E}_6$ algebra}
 \drawline(19,1)(19,9)
 \drawline(19,11)(19,20)\put(19,21){\circle*{2}}
\multiputlist(0,0)(10,0){\circle*{2},\circle{2},\circle*{2},\circle{2},\circle*{2}}
 \put(19,10){\circle{2}}
 \drawline(0,0)(8,0)
\drawline(10,0)(18,0)
\drawline(20,0)(28,0)\drawline(30,0)(38,0)
\multiputlist(0,3)(10,0){$ \alpha_2$,$ \alpha_1$,$\quad \gamma$, $ \beta_1$,$ \beta_2$}
\multiputlist(0,-3)(10,0){$ n_1$,$ n_2$,$\quad n_0$, $ n_4$,$ n_3$}
\put(21,10){$\delta_1$}
\put(21,20){$\delta_2$}
\put(14,10){$n_6$}
\put(14,20){$n_5$}
 \end{picture}
\end{center}
The category of non-degenerate $*$-representations of $\mathcal{A}_{\widetilde{E}_6,\chi}$ is equivalent to the category of non-degenerate locally-scalar representations of the graph $\widetilde{E}_6$ with the character and generalized dimension given on the following picture:

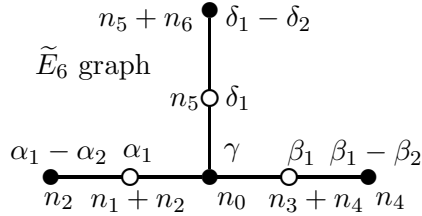
\begin{center}
 \setlength{\unitlength}{3pt}
 \begin{picture}(50,25)(-1,-3)
\thicklines
\put(-3,13){$\widetilde{E}_6$ graph}
 \drawline(19,1)(19,9)
 \drawline(19,11)(19,20)\put(19,21){\circle*{2}}
\multiputlist(0,0)(10,0){\circle*{2},\circle{2},\circle*{2},\circle{2},\circle*{2}}
 \put(19,10){\circle{2}}
 \drawline(0,0)(8,0)
\drawline(10,0)(18,0)
\drawline(20,0)(28,0)\drawline(30,0)(38,0)
\multiputlist(0,3)(10,0){$ \alpha_1-\alpha_2$,$ \alpha_1$,$\quad \gamma$, $ \beta_1$,$ \beta_1-\beta_2$}
\multiputlist(0,-3)(10,0){$ n_2$,$ n_1+ n_2$,$\quad n_0$, $ \quad n_3+n_4$,$ \quad n_4$}
\put(21,9){$\delta_1$}
\put(21,19){$\delta_1-\delta_2$}
\put(14,9){$n_5$}
\put(5,19){$n_5+n_6$}
 \end{picture}
\end{center}

It is easy to see that transition matrix  $M_f$ such that $M_f(\chi)=(x_1,x_2,\ldots, x_6,x_0)$ is block-diagonal \[ M_f= T_{6,1} \oplus T_{6,1}\oplus T_{6,1}\oplus 1,\] where \[T_{6,1} = \left(%
\begin{array}{cc}
   1 & -1 \\
  1 & 0  \\
  \end{array}%
\right).\]
The transition matrix $M_d$ such that $M_d(d_1,\ldots, d_6, d_0)^T = (n_1, \ldots, n_6, n_0)^T$ is also  block-diagonal \[ M_d^{-1} = T_{6,2}\oplus T_{6,2}\oplus T_{6,2}\oplus 1,\]  where \[ T_{6,2}= \left(%
\begin{array}{ccc}
  0 & 1  \\
  1 & 1 \\
\end{array}%
\right).
\]

For  $\Pi\in {\rep}(G,d,f)$ with sincere $d$ the Coxeter map $\overset{\bullet}{C}\overset{\circ}{C}$ transform character $x= (x_1,\ldots,x_7, x_0)$ by multiplying from the left the vector-column $x$ on the matrix $M_c= SJS^{-1}$ where the Jordan form
\[J = (-1)\oplus J_2(1) \oplus z \oplus z\oplus \overline{z}\oplus \overline{z},\] where \[J_2(1)= \left(%
\begin{array}{cc}
  1 & 1 \\
  0 & 1 \\
\end{array}%
\right),\] and $z= 1/2 (-1 - i \sqrt{3})$.

Specializing Theorem~\ref{genform} to the case  $G= \widetilde{E}_6$ we get  Theorem~\ref{ineqe6} below. In view of the symmetry of the graph $\widetilde{E}_6$ the $C$-series and  $\delta$-series of the coordinate vectors corresponding to vertices $3$ and  $5$ differ from the $C$-series and $\delta$-series of the coordinate vector corresponding to vertex $1$ by a transposition of coordinates from the set $\set{(1,4),  (1,5)}$; orbits of  the coordinate vectors corresponding to vertices $4, 6$ differ from the $C$-series and $\delta$-series of the coordinate vector corresponding to vertex $2$ by a transposition of coordinates from the set $\set{(2,4),  (2,6)}$. 
We will write $\tau v$ to denote vector obtained by permuting coordinates of vector $v$ according to permutations $\tau$. Only  the $C$-series of the roots $(1,0,0,0,0,0,0)$, $(0,1,0,0,0,0,0)$ and $(0,0,0,0,0,0,1)$ are listed in the appendix. The rest of  $C$-series can be obtained by applying corresponding permutations. The notation $u \ge_t v$ used in theorem below has been introduced  before Theorem~\ref{genform}. 

Vectors $d^{(t)}_1, d^{(t)}_2, \ldots, $ represent  roots from the $C$-series of coordinate vector $\varepsilon_t$ corresponding to the vertex $t$ of the graph listed  in the same  order they appear in the orbit of the Coxeter transformation $c$. Number $k_t$ denotes the minimal number of the Coxeter transformations $\overset{\bullet}{c}$ and $\overset{\circ}{c}$ necessary to apply (in alternating order) to get from coordinate vector $\varepsilon_t$ to a non-degenerate root.   

\begin{theorem}\label{ineqe6}
Put    $d^{(t)}_k = v^{(t)}_{k \mod m_t}+
[\frac{k}{m_t}]\delta$ for $k\ge k_t$ and $t\in \{ 0, 1, 2\}$. Here  $v^{(t)}_s$ is the
$s$-th vector in the set $C_t$  from Table~\ref{tabE6}, $m_t = |C_t|$ and $k_0 = 4$, $k_1=14$, $k_2=7$.

If the vector $\chi\not\in h_{\widetilde{E}_6}$,  i.e $\chi$ does not satisfy the equation  \begin{equation}\label{hypE6} \alpha_1+\alpha_2+\beta_1+\beta_2
+\delta_1+\delta_2 = 3 \gamma \end{equation}  then
 the algebra  $\mathcal{A}_{\widetilde{E}_6,\chi}$ has an irreducible non-degenerate
representation   in a generalized dimensions $v$ if and only if two conditions hold 
\begin{enumerate}
\item 
for some $t\in\{0,1, 2\}$ and some $k\ge
k_t$ $$v= \tau M_d d^{(t)}_k$$ where permutation $\tau\in \{ e,  (1,3),  (1,5), (2,4), (2,6)\}$ and $\tau=e$ if $t=0$.  
\item  $${\rm A}_{t,k} \tau \chi \ge_t 0$$ where the matrix ${\rm A}_{t,k}$ is
 taken from Table~\ref{tabE6}.
 \end{enumerate} The permutation $\tau$ and numbers $t$ and  $k$ are determined uniquely.    The irreducible representation of the algebra $ \mathcal{A}_{\widetilde{E}_6,\chi}$ corresponding to the above generalized dimension and $\chi$ is unique.

If $\chi$  satisfies (\ref{hypE6}) then irreducible non-degenerate
representations of $\mathcal{A}_{\widetilde{E}_6,\chi}$
may exist only in the generalized dimension $\delta_{alg}(E_6)$. They exist if
and only if $\chi$ satisfies conditions $H_{\widetilde{E}_6}$ from Table~\ref{tabE6}.    \end{theorem}

\section{ Representations of  $\mathcal{A}_{\widetilde{E}_7,\chi}$. }

The parameters $\chi=(\alpha_1, \alpha_2, \alpha_3, \beta_1, \beta_2, \beta_3, \delta, \gamma)$ of the algebra $\mathcal{A}_{\widetilde{E}_7,\chi}$ and the vector of generalized dimension $d=(d_0, d_1, \ldots, d_7)$ will be plotted on the associated graph according to the following picture:
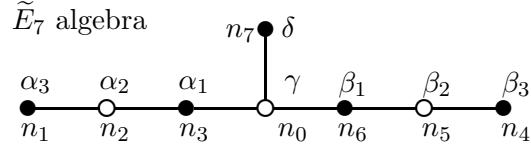
\begin{center}
\setlength{\unitlength}{3pt}
\begin{picture}(69,20)(-1,-1)
\thicklines
\put(-3,10){$\widetilde{E}_7$ algebra}
 \drawline(29,1)(29,9)
\multiputlist(0,0)(10,0){\circle*{2},\circle{2},\circle*{2},\circle{2},\circle*{2},\circle{2},\circle*{2}}
 \put(29,10){\circle*{2}}
 \drawline(0,0)(8,0)
\drawline(10,0)(18,0)
\drawline(20,0)(28,0)\drawline(30,0)(38,0)\drawline(40,0)(48,0)\drawline(50,0)(58,0)
\multiputlist(0,3)(10,0){$ {\alpha_3}$,$ \alpha_2$,$ \alpha_1$, $\quad \gamma$,$ \beta_1$, $ \beta_2$, $ \beta_3$}
\multiputlist(0,-3)(10,0){$ n_1$,$ n_2$,$ n_3$, $\quad n_0$,$ n_6$, $ n_5$, $ n_4$}
\put(31,9){$\delta$}
\put(24,9){$n_7$}
\end{picture}
\end{center}
The category of non-degenerate $*$-representations of $\mathcal{A}_{\widetilde{E}_7,\chi}$ is equivalent to the category of non-degenerate locally-scalar representations of the graph $\widetilde{E}_7$ with the character and generalized dimension given on the following picture:

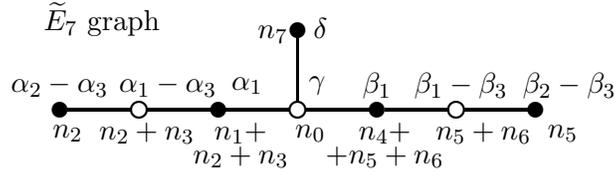
\begin{center}
\setlength{\unitlength}{3pt}
\begin{picture}(69,20)(-1,-6)
\thicklines
\put(-3,10){$\widetilde{E}_7$ graph}
 \drawline(29,1)(29,9)
\multiputlist(0,0)(10,0){\circle*{2},\circle{2},\circle*{2},\circle{2},\circle*{2},\circle{2},\circle*{2}}
 \put(29,10){\circle*{2}}
 \drawline(0,0)(8,0)
\drawline(10,0)(18,0)
\drawline(20,0)(28,0)\drawline(30,0)(38,0)\drawline(40,0)(48,0)\drawline(50,0)(58,0)
\multiputlist(-1,3)(10,0){\small ${\alpha_2-{\alpha_3}}$,\small $\quad\quad   \alpha_1-\alpha_3$,\small $ \quad \quad \alpha_1$, \small $ \quad \gamma$,\small $  \beta_1$, \small $ \beta_1-\beta_3$, \small $\quad  \quad \beta_2-\beta_3$}
\multiputlist(0,-3)(10,0){\small $n_2$,\small $n_2+n_3$,\small \quad $n_1+$, \small $n_0$,\small $ n_4+$, \small \quad $ n_5+n_6$, \small \quad $ n_5$}
\multiputlist(0,-6)(10,0){ {},{},\small \quad $ n_2+n_3$, {},\small $+n_5+n_6$}
\put(31,9){\small $\delta$}
\put(24,9){\small $n_7$}
\end{picture}
\end{center}

It is easy to see that transition matrix  $M_f$ such that $M_f(\chi)=(x_1,x_2,\ldots, x_7,x_0)$ is block-diagonal \[ M_f= \diag{(T_1, T_1, 1, 1)},\] where \[T_1 = \left(%
\begin{array}{ccc}
  0 & 1 & -1 \\
  1 & 0 & -1 \\
  1 & 0 & 0 \\
\end{array}%
\right).\]
The transition matrix $M_d$ such that $M_d(d_0,d_1,\ldots, d_7)^T = (n_0, n_1, \ldots, n_7)^T$ is also  block-diagonal \[ M_d^{-1} = \diag{(T_2, T_2, 1)},\]  where \[ T_2= \left(%
\begin{array}{ccc}
  0 & 1 & 0 \\
  0 & 1 & 1 \\
  1 & 1 & 1 \\
\end{array}%
\right).
\]

For  $\Pi\in {\rep}(G,d,f)$ with sincere $d$ the Coxeter map $\overset{\bullet}{C}\overset{\circ}{C}$ transform character $x= (x_1,\ldots,x_7, x_0)$ by multiplying from the left the vector-column $x$ on the matrix $M_c= SJS^{-1}$ where the Jordan form
\[J = \diag{(-1,-1,J_2(1),-i,i, 1/2 (-1 - i \sqrt{3}),1/2 (-1 + i \sqrt{3}))},\] where \[J_2(1)= \left(%
\begin{array}{cc}
  1 & 1 \\
  0 & 1 \\
\end{array}%
\right).\]

Specializing Theorem~\ref{genform} to the case  $G= \widetilde{E}_7$ we get  Theorem~\ref{ineqe7} below. In view of the symmetry of the graph $\widetilde{E}_7$ the $C$-series and  $\delta$-series of the coordinate vector corresponding to vertex $4$ differ from the $C$-series and $\delta$-series of the coordinate vector corresponding to vertex $1$ by the  transposition of coordinates  $\tau_1=(1,4)$; the orbit of  the coordinate vector corresponding to $5$ differ from the $C$-series and $\delta$-series of the coordinate vector corresponding to vertex $2$ by a transposition of coordinates  $\tau_2=(2,5)$; the orbit of  the coordinate vector corresponding to $6$ differs from the $C$-series and $\delta$-series of the coordinate vector corresponding to vertex $3$ by a transposition of coordinates  $\tau_3=(3,6)$. 
We will write $\tau v$ to denote vector obtained by permuting coordinates of vector $v$ according to permutations $\tau$. Only  the $C$-series of the roots $(1,0,0,0,0,0,0,0)$, $(0,1,0,0,0,0,0,0)$, $(0,0,1,0,0,0,0,0)$, $(0,0,0,0,0,0,0,1)$ are listed in the appendix. The rest of  $C$-series can be obtained by applying corresponding permutations. The notation $u \ge_t v$ used in theorem below has been introduced  before Theorem~\ref{genform}.  

Vectors $d^{(t)}_1, d^{(t)}_2, \ldots $ represent  roots from the $C$-series of coordinate vector $\varepsilon_t$ corresponding to the vertex $t$ of the graph listed  in the same  order they appear in the orbit of the Coxeter transformation $c$. Number $k_t$ denotes the minimal number of the Coxeter transformations $\overset{\bullet}{c}$ and $\overset{\circ}{c}$ necessary to apply (in alternating order) to get from coordinate vector $\varepsilon_t$ to a non-degenerate root.

 \begin{theorem}\label{ineqe7}
Put  $d^{(t)}_k = v^{(t)}_{k \mod m_t}+
 \epsilon_t [\frac{k}{m_t}]\delta$ for $k\ge k_t$ and $t\in \{ 0, 1, 2, 3, 7, 8 \}$  where
 $v^{(t)}_s$ is the $s$-th vector in the set $C_t$
 from Table~\ref{tabE7}, $m_t = |C_t|$ and $k_0 = 4$, $k_1=27$, $k_2=14$,  $k_3=9$,  $k_7=11$,  $k_8=6$,
 $\epsilon_j=1$ if $j\not\in\{2,5\}$ and  $\epsilon_2=\epsilon_5=2$.

If the vector $\chi\not\in h_{\widetilde{E}_7}$, i.e $\chi$ does not satisfy the equation 
\begin{equation}\label{hypE7} \alpha_1+\alpha_2+\alpha_3+
\beta_1+\beta_2 +\beta_3+ 2\delta_1 = 4 \gamma
\end{equation}  then
 the algebra  $\mathcal{A}_{\widetilde{E}_7,\chi}$ has an irreducible non-degenerate
representation   in a generalized dimensions $v$ if and only if two conditions hold 

\begin{enumerate}
\item 
for some $t\in\{0,1, 2\}$ and some $k\ge
k_t$ $$v= \tau M_d d^{(t)}_k$$ where transposition $\tau\in \{ e, (1,4),  (2,5),  (3,6)\}$ and $\tau=e$ if $t=0$,
\item  $${\rm A}_{t,k} \tau \chi \ge_t 0$$
 where the matrix ${\rm A}_{t,k}$ is taken from Table~\ref{tabE7}. 
 \end{enumerate}
 
 The permutation $\tau$ and numbers $t$ and  $k$ are determined uniquely.    The irreducible representation of the algebra $ \mathcal{A}_{\widetilde{E}_7,\chi}$ corresponding to the above generalized dimension and $\chi$ is unique.

If $\chi$  satisfies (\ref{hypE7}) then irreducible non-degenerate representations of
$\mathcal{A}_{\widetilde{E}_7,\chi}$ may exist only in the generalized
dimension $\delta_{alg}(\widetilde{E}_7)$. They exist if and only
if $\chi$ satisfies conditions $H_{\widetilde{E}_7}$ from Table~\ref{tabE7}.    \end{theorem}

\section{ Representations of  $\mathcal{A}_{\widetilde{E}_8,\chi}$. }

The parameters $\chi=(\alpha_1, \alpha_2, \alpha_3,  \alpha_4, \alpha_5, \beta_1, \beta_2,  \delta, \gamma)$ of the algebra $\mathcal{A}_{\widetilde{E}_8,\chi}$ and the vector of generalized dimension $d=(d_1, \ldots, d_8, d_0)$ will be plotted on the associated graph according to the following picture:
\begin{center}
\setlength{\unitlength}{4pt}
\begin{picture}(62,25)(-1,-1)
\thicklines
\put(-3,10){\small $\widetilde{E}_8$ algebra }
 \drawline(49,1)(49,9)
\multiputlist(0,0)(10,0){\circle*{2},\circle{2},\circle*{2},\circle{2},\circle*{2},\circle{2},\circle*{2},\circle{2}}
 \put(49,10){\circle*{2}}
 \drawline(0,0)(8,0)
\drawline(10,0)(18,0)
\drawline(20,0)(28,0)\drawline(30,0)(38,0)\drawline(40,0)(48,0)\drawline(50,0)(58,0)\drawline(60,0)(68,0)
\multiputlist(-1,3)(10,0){\small ${\alpha_5}$,\small $\alpha_4$,\small $  \alpha_3$, \small $ \alpha_2$, \small $ \alpha_1$,\small $ \quad \gamma$,\small $  \beta_1$, \small $ \beta_2$}
\multiputlist(-1,-3)(10,0){\small $n_1$,\small $n_2$,\small   $n_3$, \small  $ n_4$,\small   $ n_5$, \small \quad $ n_0$, \small  $n_7$, \small  $n_6$}
\put(51,9){\small $\delta$}
\put(45,9){\small $n_8$}
\end{picture}
\end{center}
The category of non-degenerate $*$-representations of $\mathcal{A}_{\widetilde{E}_8,\chi}$ is equivalent to the category of non-degenerate locally-scalar representations of the graph $\widetilde{E}_8$ with the character and generalized dimension given on the following picture:

\begin{center}
\setlength{\unitlength}{4pt}
\begin{picture}(62,20)(1,-9)
\thicklines
\put(-3,10){$\widetilde{E}_8$ graph}
 \drawline(49,1)(49,9)
\multiputlist(0,0)(10,0){\circle*{2},\circle{2},\circle*{2},\circle{2},\circle*{2},\circle{2},\circle*{2},\circle{2}}
 \put(49,10){\circle*{2}}
 \drawline(0,0)(8,0)
\drawline(10,0)(18,0)
\drawline(20,0)(28,0)\drawline(30,0)(38,0)\drawline(40,0)(48,0)\drawline(50,0)(58,0)\drawline(60,0)(68,0)
\multiputlist(-1,3)(10,0){\small ${\alpha_3-{\alpha_4}}$,\small $\alpha_2-\alpha_4$,\small $  \alpha_2-\alpha_5$, \small $ \alpha_1-\alpha_5$, \small $ \alpha_1$,\small $ \quad \gamma$,\small $  \beta_1$, \small $ \beta_1-\beta_2$}
\multiputlist(-1,-3)(10,0){\small $n_3$,\small $n_3+n_4$,\small  \quad $n_2+n_3$, \small \quad $ n_2+n_3$,\small \quad\quad  $ n_1+n_2$, \small \quad $ n_0$, \small \quad $n_6+n_7$, \small \quad $n_7$}
\multiputlist(0,-6)(10,0){ {},{},\small  $+n_4$, \small {$+n_4+$},\small $+n_3+$}
\multiputlist(0,-9)(10,0){ {},{},{}, \small \quad{$n_5$},\small \quad \quad $n_4+n_5$}
\put(51,9){\small $\delta$}
\put(45,9){\small $n_8$}
\end{picture}
\end{center}

It is easy to see that transition matrix  $M_f$ such that $M_f(\chi)=(x_1,x_2,\ldots, x_8,x_0)$ is block-diagonal \[ M_f= \diag{(T_2, T_3, 1, 1)},\] where \[T_2 = \left(%
\begin{array}{ccccc}
  0 & 0 & 1  & -1 & 0 \\
  0& 1 & 0 & -1 & 0\\
  0 & 1 & 0 & 0 & -1 \\
  1 & 0 & 0 & 0 & -1\\
  1 & 0 & 0 & 0 & 0\\
\end{array}%
\right),
T_3 = \left(%
\begin{array}{cc}
  1 & -1  \\
  1 & 0 \\
  \end{array}%
\right).    \]
The transition matrix $M_d$ such that $M_d(d_1,\ldots, d_8, d_0)^T = (n_1, \ldots, n_8, n_0)^T$ is also  block-diagonal \[ M_d^{-1} = \diag{(T_4, T_5, 1, 1)},\]  where \[ T_4= \left(%
\begin{array}{ccccc}
  0 & 0 & 1 & 0 & 0\\
  0 & 0 & 1 & 1 & 0 \\
  0 & 1 & 1 & 1 & 0 \\
  0 & 1 & 1 & 1 & 1 \\
  1 & 1 & 1 & 1 & 1 \\
\end{array}%
\right),
T_5 = \left(%
\begin{array}{cc}
  0 & 1  \\
  1 & 1 \\
  \end{array}%
\right).
\]

For  $\Pi\in {\rep}(G,d,f)$ with sincere $d$ the Coxeter map $\overset{\bullet}{C}\overset{\circ}{C}$ transform character $x= (x_1,\ldots,x_7, x_0)$ by multiplying from the left the vector-column $x$ on the matrix $M_c= SJS^{-1}$ where the Jordan form
\[J = \diag{(-1,J_2(1), 1/2 (-1 - i \sqrt{3}),1/2 (-1 + i \sqrt{3}), \zeta, \zeta^2, \zeta^3, \zeta^4)},\] where \[J_2(1)= \left(%
\begin{array}{cc}
  1 & 1 \\
  0 & 1 \\
\end{array}%
\right),\]
and $\zeta$ is a prime root of unity of degree $5$

Specializing Theorem~\ref{genform} to the case  $G= \widetilde{E}_8$ we get  Theorem~\ref{ineqe8} below. 
 The notation $u \ge_t v$ used in theorem below has been introduced  before Theorem~\ref{genform}.
 
  Vectors $d^{(t)}_1, d^{(t)}_2, \ldots $ represent  roots from the $C$-series of coordinate vector $\varepsilon_t$ corresponding to the vertex $t$ of the graph listed  in the same  order they appear in the orbit of the Coxeter transformation $c$. Number $k_t$ denotes the minimal number of the Coxeter transformations $\overset{\bullet}{c}$ and $\overset{\circ}{c}$ necessary to apply (in alternating order) to get from coordinate vector $\varepsilon_t$ to a non-degenerate root, i.e.  all vectors  $d^{(t)}_k$ are non-degenerate for $k\ge k_t$  and degenerate otherwise.

 \begin{theorem}\label{ineqe8}
Put $d^{(t)}_k = v^{(t)}_{k \mod m_t}+
 \epsilon_t [\frac{k}{m_t}]\delta$ for $k\ge k_t$ and $t\in \{ 0, 1,\ldots, 8 \}$.
 Here  $v^{(t)}_s$ is the $s$-th vector in the set $C_t$
 from Table~\ref{tabE8}, $m_t = |C_t|$ and $(k_0, \ldots, k_8) = (10, 65, 34, 23, 14, 13, 28, 17, 19)$, $\epsilon_j=1$ if $j\not\in\{4,7\}$ and  $\epsilon_4=\epsilon_7=2$.

If the vector $\chi\not\in h_{\widetilde{E}_8}$, i.e. $\chi$ does not satisfy the equation  \begin{equation}\label{hypE8} \alpha_1+\alpha_2+\alpha_3+\alpha_4+ 2(\beta_1+\beta_2)  +3\delta_1 = 6 \gamma
\end{equation}  then
 the algebra  $\mathcal{A}_{\widetilde{E}_8,\chi}$ has an irreducible non-degenerate
representation   in a generalized dimensions $v$ if and only if two conditions hold 

\begin{enumerate}
\item for some $t\in\{0, 1,\ldots, 8 \}$ and
some $k\ge k_t$ $$v=  M_d d^{(t)}_k,$$ 
\item 
 ${\rm A}_{t,k}  \chi \ge_t 0$ where the matrix ${\rm A}_{t,k}$ is taken from
Table~\ref{tabE8}.
\end{enumerate} 

 The permutation $\tau$ and numbers $t$ and  $k$ are determined uniquely.    The irreducible representation of the algebra $ \mathcal{A}_{\widetilde{E}_8,\chi}$ corresponding to the above generalized dimension and $\chi$ is unique.
 
If $\chi$  satisfies (\ref{hypE8}) then irreducible non-degenerate representations of
$\mathcal{A}_{\widetilde{E}_8,\chi}$ may exist only in the generalized dimension
$\delta_{alg}(\widetilde{E}_8)$. They exist if and only if $\chi$ satisfies
conditions $H_{\widetilde{E}_8}$ from Table~\ref{tabE8}.    \end{theorem}

\section*{Appendix: Tables. }
Below  we present  tables 7-10 for extended Dynkin graphs $\widetilde{D}_4$, $\widetilde{E}_6$, $\widetilde{E}_7$, $\widetilde{E}_8$ correspondingly. 

For easy referencing we recall some notations from the paper used in the tables.
  
The set of roots $\Delta_G$ of graph $G$ is a disjoint union 
$\Delta_{reg}\cup \Delta_{sing}$ of the sets of regular and singular roots. Each of the sets 
$\Delta_{reg}$ and $\Delta_{sing}$  are union of elements of  $\delta$-series (cosets in  $\Delta / \delta \mathbb{Z}$). Here $\delta=\delta_G$ is a minimal imaginary root. 
The set of representatives $\Delta_f$ of all $\delta$-series is a finite set. We present in tables the  sets $ \Delta_{reg}\cap \Delta_f$. 

The set $\Delta_{sing}$  is a union of finite number of $C$-series (orbits under the action of the group generated by Coxeter transformations $\overset{\bullet}{c}$ and $\overset{\circ}{c}$) let us denote these $C$-series by  $\mathcal{C}_1, \ldots, \mathcal{C}_{m}$ (each $\mathcal{C}_j$ corresponds to $j$-vertex of the graph). We present in the tables the finite sets $C_j$ of representatives of $\delta$-series comprising $\mathcal{C}_j$. The reader should bear in mind that the set $C_j$ is a ``part'' of the orbit of coordinate vector $\varepsilon_j$ rather then the intersection $\mathcal{C}_j\cap \Delta_f$.     

 Matrices $A_{j,k}$ define linear inequalities on the coordinates of a vector $\chi$  which are necessary and sufficient conditions for $*$-algebra $A_{G,\chi}$ to have a representation in generalized dimension $d_{k}^{(j)}$ (see Theorems~\ref{ineqeD4}-\ref{ineqe8}).    

In the following tables coordinates of the root vectors $(v_1, \ldots, v_n, v_0)$
correspond  to the enumeration of vertices shown on the pictures in sections 3-6.
We will omit the parentheses and commas in the vectors to shorten notations.
If a coordinate is not a decimal digit  it will be put in parentheses.
{\small
\section{Table. $C$-orbits and  matrices ${\rm A}_{i,j}$ for $\widetilde{D}_4$. }\label{tabD4}

\noindent 1.  $\Delta$ consists of $25$ $\delta$-orbits.

\noindent 2.  $\Delta$ is a union of 5 $C$-series consisting of singular roots and   7 $\delta$-series of regular roots. $ \Delta_f\cap \Delta_{reg}=$  \begin{eqnarray*} &\pm&\{  1 0 0 1 1, 1 0 1 0 1, 1 1 0 0 1, 00111, \\ & & 0 1 0 1 1, 0 1 1 01  \}. \end{eqnarray*} 

\noindent 3.  $C_1:  10000, 10001, 01111, 01112$.
$\overset{\circ}{C} v_0 =v_1, \ldots, \overset{\bullet}{C} v_{3} =v_0 + \delta$.

\noindent 4. $C_0: 00001, 11111$.  $\overset{\bullet}{C} v_0 =v_1,  \overset{\circ}{C} v_{1} =v_0 + \delta$.

\noindent 5.

$$ \begin{array}{cc}
  \rm{A}_{1,k}= \qquad\qquad\qquad\qquad\qquad\qquad\qquad & \rm{A}_{0,k}= \qquad\qquad\qquad\qquad\qquad\qquad \\
   { D_{1, 5}\begin{cases}  C^{s-3}\overset{\bullet}{C} M_f  &
\text{if $k=2s$}\\
C^{s-3} M_f  & \text{if $k=2s-1$}
\end{cases}},  & {D_{0,2}\begin{cases}  C^{s-1} M_f  &
\text{if $k=2s$}\\
C^{s-1} \overset{\bullet}{C}M_f  & \text{if $k=2s+1$}
\end{cases}};
 \\
\end{array}
 $$

$$ \begin{array}{cc}
  D_{1, 5} = \qquad\qquad\qquad\qquad\qquad\qquad\qquad & D_{0,2}= \qquad\qquad\qquad\qquad\qquad\qquad \\
   \tiny{\left(%
\begin{array}{ccccc}
 -2 & -1 & -1 & -1 & 3 \\
 -1 & 0 & -1 & -1 & 2 \\
 -1 & -1 & 0 & -1 & 2 \\
 -1 & -1 & -1 & 0 & 2 \\
 -3 & -1 & -1 & -1 & 4\\
\end{array}%
\right)},  & \tiny{\left(%
\begin{array}{ccccc}
 -1 & 0 & 0 & 0 & 1 \\
 0 & -1 & 0 & 0 & 1 \\
 0 & 0 & -1 & 0 & 1 \\
 0 & 0 & 0 & -1 & 1 \\
 -1 & -1 & -1 & -1 & 3\\
\end{array}%
\right)};
 \\
\end{array}
 $$

 \noindent 6.
 $\delta_{alg}(\widetilde{D}_4)= (1,1,1,1;2)$.

Hyperplane conditions $H_{\widetilde{D}_4}$:
$$
\alpha+\beta+\xi> \gamma, \alpha+\xi+\delta> \gamma, \alpha+\beta+\delta> \gamma, \beta+\xi+\delta> \gamma.
$$

\section{Table. $C$-orbits and  matrices ${\rm A}_{i,j}$  for $\widetilde{E}_6$. }\label{tabE6}

1.  $\Delta$ consists of $73$ $\delta$-orbits.

\noindent 2.  $\Delta$ is a union of 7 $C$-series consisting of singular roots and   15
 $\delta$-series of regular roots. $ \Delta_f\cap \Delta_{reg}=$ \begin{eqnarray*}   &\pm&\{ 0001111, 0011011, 0100111, 0101012, 0111001,\\ & &0111122, 0112112  \}. \end{eqnarray*}.

\noindent 3. $C_1:$ \begin{eqnarray*}   &&1000000, 1100000, 0100001, 0001011, 0011111, 0111111,\\ &&1101012, 1201012,
1211112, 1112122, 0112123, 0212123.
 \end{eqnarray*}
  $\overset{\circ}{C} v_0 =v_1, \ldots, \overset{\bullet}{C} v_{11} =v_0 + \delta$.

\noindent 4.  $C_2:$ \begin{gather*} 0100000, 1100001, 1101011, 0111112, 0112122, 1112123.
\end{gather*} $\overset{\bullet}{C} v_0 =v_1, \ldots, \overset{\circ}{C} v_{5} =v_0 + \delta$.

\noindent 5.  $C_0:$  $$0000001, 0101011, 1111112, 1212122.$$

$\overset{\circ}{C} v_0 =v_1, \ldots, \overset{\bullet}{C} v_{3} =v_0 + \delta$

\noindent 6.

$$ \begin{array}{cc}
  \rm{A}_{1,k}= \qquad\qquad\qquad\qquad\qquad\qquad\qquad & \rm{A}_{2,k}= \qquad\qquad\qquad\qquad\qquad\qquad \\
  { D_{1, 14}\begin{cases}  C^{s-7} M_f  &
\text{if $k=2s$}\\
C^{s-7}\overset{\bullet}{C} M_f  & \text{if $k=2s+1$}
\end{cases}},  & { D_{2,7}\begin{cases}  C^{4-s}\overset{\circ}{C} M_f  &
\text{if $k=2s$,}\\
C^{4-s} M_f  & \text{if $k=2s-1$.}
\end{cases}};
 \\
\end{array}
 $$

$$\rm{A}_{0,k}=D_{0,4}\begin{cases}  C^{2-s} M_f  &
\text{if $k=2s$}\\
C^{2-s}\overset{\circ}{C} M_f  & \text{if $k=2s+1$}
\end{cases}
$$

$$ \begin{array}{cc}
  D_{1, 14} = \qquad\qquad\qquad\qquad\qquad\qquad\qquad & D_{2,7}=  \qquad\qquad\qquad\qquad\qquad\qquad \\
   \tiny{ \left(%
\begin{array}{ccccccc}
 1 & -3 & 1 & -2 & 1 & -2 & 4 \\
 2 & -5 & 2 & -3 & 2 & -3 & 6 \\
 1 & -2 & 0 & -1 & 1 & -2 & 3 \\
 1 & -4 & 1 & -3 & 2 & -3 & 6 \\
 1 & -2 & 1 & -2 & 0 & -1 & 3 \\
 1 & -4 & 2 & -3 & 1 & -3 & 6 \\
 2 & -7 & 2 & -4 & 2 & -4 & 9 \\
\end{array}%
\right)},  & \tiny{ \left(%
\begin{array}{ccccccc}
 1 & -1 & 0 & -1 & 0 & -1 & 2 \\
 2 & -3 & 1 & -2 & 1 & -2 & 4 \\
 1 & -1 & 0 & 0 & 1 & -1 & 1 \\
 1 & -2 & 1 & -1 & 1 & -2 & 3 \\
 1 & -1 & 1 & -1 & 0 & 0 & 1 \\
 1 & -2 & 1 & -2 & 1 & -1 & 3 \\
 2 & -4 & 1 & -2 & 1 & -2 & 5\\
\end{array}%
\right)};
 \\
\end{array}
 $$

$$ D_{0,4}= \tiny{\left(%
\begin{array}{ccccccc}
  0 & -1 & 0 & 0 & 0 & 0 & 1 \\
 0 & -1 & 1 & -1 & 1 & -1 & 2 \\
 0 & 0 & 0 & -1 & 0 & 0 & 1 \\
 1 & -1 & 0 & -1 & 1 & -1 & 2 \\
 0 & 0 & 0 & 0 & 0 & -1 & 1 \\
 1 & -1 & 1 & -1 & 0 & -1 & 2 \\
 1 & -2 & 1 & -2 & 1 & -2 & 4 \\
\end{array}%
\right)} $$

\noindent 7.  $\delta_{alg}(\widetilde{E}_6)= (1,1;1,1;1,1;3)$.

Hyperplane conditions $H_{\widetilde{E}_6}$:

\begin{gather*}\alpha_1 + \beta_1 > \gamma, \alpha_1 + \beta_2 + \delta_2 > \gamma, \alpha_1 + \delta_1 > \gamma, \alpha_2 + \beta_1 + \delta_2 > \gamma, \\
 \alpha_2 + \beta_2 + \delta_1 > \gamma, \beta_1 + \delta_1 > \gamma, \alpha_1 + \alpha_2 + \beta_1 + \beta_2 + \delta_2 >
  2 \gamma, \\
   \alpha_1 + \alpha_2 + \beta_1 + \delta_1 >
  2 \gamma, \alpha_1 + \alpha_2 + \beta_2 + \delta_1 + \delta_2 > 2 \gamma, \alpha_1 + \beta_1 + \beta_2 + \delta_1 >
   2 \gamma, \\  \alpha_1 + \beta_1 + \delta_1 + \delta_2 >
  2 \gamma, \alpha_2 + \beta_1 + \beta_2 + \delta_1 + \delta_2 > 2 \gamma.
  \end{gather*}

\section{Table. $C$-orbits and  matrices ${\rm A}_{i,j}$ for $\widetilde{E}_7$. }\label{tabE7}

1.  $\Delta$ consists of $63$ $\delta$-orbits.

\noindent 2. $\Delta$ is a union of 7 $C$-series consisting of singular roots and   21 $\delta$-series of regular roots. $ \Delta_f\cap \Delta_{reg}=$
 \begin{align*} \pm \{  00011111, 00100112, 00101101, 00112212, 01100101,\\
 01101111, 01111212,
01201223, 01211112, 01212313  \}. \end{align*}  

\noindent 3. $C_1:$ \begin{eqnarray*}   &&10000000, 11000000, 01100000, 00100001, 00000111,\\
 &&00001111, 00111101,
01111101, 11101111, 11100112,\\
 &&01200112, 01201112, 11111212, 11112212,
01212212,\\
 &&01211213, 11201223, 12201223, 12311213, 12312213,\\
 &&12212323,
11212324, 01312324, 02312324.
\end{eqnarray*}
 $\overset{\circ}{C} v_0 =v_1, \ldots, \overset{\bullet}{C} v_{23} =v_0 + \delta$.

\noindent 4. $C_2:$  \begin{eqnarray*}  &&01000000, 11100000, 11100001, 01100111, 00101112,\\
 &&00111212, 01112212,
11212212, 12211213, 12301223,\\
 &&12301224, 12311324, 12313324, 12323424,
12323425,\\
 &&12413435, 13412436, 23512436, 24513436, 24523536,\\
 &&23524537,
13524647, 13524648, 23624648.
\end{eqnarray*}
 $\overset{\bullet}{C} v_0 =v_1, \ldots, \overset{\circ}{C} v_{23} =v_0 +2 \delta$.

\noindent 5. $C_3:$ \begin{eqnarray*}    &&00100000, 01100001, 11100111, 11101112, 01211212,\\
 &&01212213, 11212323, 12212324.
\end{eqnarray*}
  $\overset{\circ}{C} v_0 =v_1, \ldots, \overset{\bullet}{C} v_{7} =v_0 + \delta$.

\noindent 6. $C_7:$ \begin{eqnarray*}    &&00000010, 00000011, 00100101, 01101101, 11111111,\\
 &&11111112, 01201212,
01201213, 11211223, 12212223,\\
 &&12312313, 12312314.
\end{eqnarray*}

$\overset{\bullet}{C} v_0 =v_1, \ldots, \overset{\circ}{C} v_{11} =v_0 + \delta$.

\noindent 7. $C_0:$ \begin{eqnarray*}    &&00000001, 00100111, 01101112, 11211212, 12212213,\\ &&12312323.
  \end{eqnarray*}
    $\overset{\bullet}{C} v_0 =v_1, \ldots, \overset{\circ}{C} v_{5} =v_0 + \delta$.

\noindent 8.
$$ \begin{array}{cc}
  \rm{A}_{1,k}= \qquad\qquad\qquad\qquad\qquad\qquad\qquad & \rm{A}_{2,k}= \qquad\qquad\qquad\qquad\qquad\qquad \\
   { D_{1, 27}\begin{cases}  C^{s-14} \overset{\bullet}{C} M_f  &
\text{if $k=2s$}\\
C^{s-14} M_f  & \text{if $k=2s-1$}
\end{cases} },  & { D_{2, 14}\begin{cases}  C^{7-s}  M_f  &
\text{if $k=2s$}\\
\overset{\bullet}{C} C^{s-7} M_f  & \text{if $k=2s+1$}
\end{cases} };
 \\
\end{array}
 $$

$$ \begin{array}{cc}
  \rm{A}_{3,k}= \qquad\qquad\qquad\qquad\qquad\qquad\qquad & \rm{A}_{7,k}= \qquad\qquad\qquad\qquad\qquad\qquad \\
   { D_{3, 9}\begin{cases}  C^{s-5}\overset{\bullet}{C}  M_f  &
\text{if $k=2s$}\\
 C^{s-5} M_f  & \text{if $k=2s-1$}
\end{cases} },  & { D_{7, 11}\begin{cases}  \overset{\circ}{C}C^{s-6}  M_f  &
\text{if $k=2s$}\\
 C^{6-s} M_f  & \text{if $k=2s-1$}
\end{cases} };
 \\
\end{array}
 $$

$$\rm{A}_{8,k}=
D_{8, 6}\begin{cases}  C^{s-3}  M_f  &
\text{if $k=2s$}\\
 C^{s-3} \overset{\bullet}{C} M_f  & \text{if $k=2s+1$}
\end{cases}
$$

$$ \begin{array}{cc}
  D_{1,27} = \qquad\qquad\qquad\qquad\qquad\qquad\qquad & D_{2,14} = \qquad\qquad\qquad\qquad\qquad\qquad \\
   \tiny{ \left(%
\begin{array}{cccccccc}
 -1 & 2 & -4 & -1 & 2 & -3 & -2 & 5 \\
 -2 & 4 & -7 & -2 & 4 & -5 & -3 & 8 \\
 -2 & 5 & -10 & -2 & 5 & -7 & -5 & 12 \\
 -1 & 2 & -3 & -1 & 1 & -2 & -2 & 4 \\
 -1 & 3 & -6 & -2 & 3 & -5 & -3 & 8 \\
 -2 & 4 & -9 & -3 & 5 & -7 & -5 & 12 \\
 -2 & 3 & -6 & -1 & 3 & -5 & -3 & 8 \\
 -3 & 6 & -13 & -3 & 6 & -9 & -6 & 16 \\
\end{array}%
\right) },  & \tiny{ \left(%
\begin{array}{cccccccc}
 -1 & 1 & -1 & -1 & 2 & -2 & -1 & 2 \\
 -1 & 2 & -3 & -2 & 3 & -4 & -2 & 5 \\
 -1 & 3 & -4 & -3 & 4 & -5 & -2 & 6 \\
 0 & 0 & -1 & -1 & 1 & -1 & -1 & 2 \\
 0 & 1 & -2 & -1 & 2 & -3 & -2 & 4 \\
 -1 & 2 & -3 & -2 & 3 & -5 & -2 & 6 \\
 -1 & 2 & -3 & -1 & 2 & -3 & -1 & 4 \\
 -2 & 4 & -5 & -3 & 4 & -6 & -3 & 8 \\
\end{array}%
\right) };
 \\
\end{array}
 $$

$$ \begin{array}{cc}
  D_{3,9} = \qquad\qquad\qquad\qquad\qquad\qquad\qquad & D_{7,11} = \qquad\qquad\qquad\qquad\qquad\qquad \\
   \tiny{ \left(%
\begin{array}{cccccccc}
  1 & -1 & 1 & 1 & -1 & 2 & 1 & -2 \\
 1 & -2 & 3 & 2 & -3 & 4 & 2 & -4 \\
 2 & -3 & 5 & 2 & -4 & 6 & 4 & -7 \\
 0 & 0 & 1 & 0 & -1 & 2 & 1 & -2 \\
 1 & -1 & 3 & 1 & -2 & 3 & 2 & -4 \\
 2 & -3 & 5 & 1 & -3 & 5 & 3 & -6 \\
 1 & -2 & 3 & 1 & -2 & 3 & 3 & -4 \\
 3 & -4 & 6 & 2 & -5 & 7 & 5 & -8 \\
\end{array}%
\right)},  & \tiny{ \left(%
\begin{array}{cccccccc}
  0 & 1 & -1 & -1 & 1 & -1 & 0 & 1 \\
 -1 & 2 & -2 & -1 & 1 & -2 & -1 & 3 \\
 -1 & 2 & -3 & -1 & 2 & -4 & -1 & 5 \\
 -1 & 1 & -1 & 0 & 1 & -1 & 0 & 1 \\
 -1 & 1 & -2 & -1 & 2 & -2 & -1 & 3 \\
 -1 & 2 & -4 & -1 & 2 & -3 & -1 & 5 \\
 -1 & 2 & -3 & -1 & 2 & -3 & -1 & 4 \\
 -1 & 3 & -5 & -1 & 3 & -5 & -2 & 7 \\
\end{array}%
\right)};
 \\
\end{array}
 $$

 $$D_{8,6} = \tiny{\left(%
\begin{array}{cccccccc}
  0 & 0 & -1 & 0 & 0 & 0 & 0 & 1 \\
 0 & 0 & -1 & 0 & 1 & -1 & -1 & 2 \\
 0 & 1 & -2 & -1 & 2 & -2 & -1 & 3 \\
 0 & 0 & 0 & 0 & 0 & -1 & 0 & 1 \\
 0 & 1 & -1 & 0 & 0 & -1 & -1 & 2 \\
 -1 & 2 & -2 & 0 & 1 & -2 & -1 & 3 \\
 -1 & 1 & -1 & -1 & 1 & -1 & -1 & 2 \\
 -1 & 2 & -3 & -1 & 2 & -3 & -2 & 5 \\
\end{array}%
\right)}
 $$
\item
 $\delta_{alg}= (1, 1, 1, 1, 1, 1, 2; 4)$.

Hyperplane conditions $H_{\widetilde{E}_7}$:
\begin{eqnarray*} &&\alpha_1 + \beta_2 > \gamma, \delta + \alpha_1 > \gamma, \alpha_2 + \beta_1 > \gamma,
 \delta + \alpha_2 + \beta_3 > \gamma,
 \\ &&\delta + \alpha_3 + \beta_2 > \gamma,   \delta + \beta_1 > \gamma,
 \delta + \alpha_1 + \alpha_2 + \beta_2 >
  2 \gamma,
  \\ &&\delta + \alpha_1 + \alpha_3 + \beta_1 >
  2 \gamma,  \delta + \alpha_1 + \alpha_3 + \beta_2 + \beta_3 > 2 \gamma, \delta + \alpha_1 + \beta_1 + \beta_3 >
  2 \gamma,
  \\  &&\delta + \alpha_2 + \alpha_3 + \beta_1 + \beta_3 > 2 \gamma, \delta + \alpha_2 + \beta_1 + \beta_2 >
  2 \gamma,
  \\ &&\delta + \alpha_1 + \alpha_2 + \alpha_3 + \beta_1 + \beta_2 > 3 \gamma,
2 \delta + \alpha_1 + \alpha_2 + \alpha_3 + \beta_2 + \beta_3 > 3 \gamma,
\\ &&\delta + \alpha_1 + \alpha_2 + \beta_1 + \beta_2 + \beta_3 > 3 \gamma, 2 \delta + \alpha_1 + \alpha_2 + \beta_1 + \beta_3 >   3 \gamma, \\
&&2 \delta + \alpha_1 + \alpha_3 + \beta_1 + \beta_2 >
  3 \gamma,
2 \delta + \alpha_2 + \alpha_3 + \beta_1 + \beta_2 + \beta_3 > 3 \gamma
  \end{eqnarray*}

\section{Table. $C$-orbits and  matrices ${\rm A}_{i,j}$ for $\widetilde{E}_8$.} \label{tabE8}

\noindent 1. $\Delta$ consists of $241$ $\delta$-orbits.

\noindent 2. $\Delta$ is a union of 8 $C$-series consisting of singular roots and   29 $\delta$-series of
regular roots. $ \Delta_f\cap \Delta_{reg}=$ \begin{eqnarray*} &&\pm \{  000010112, 000111111, 000121223, 001111101,\\
 &&001121213, 001220112,
001232324, 011110101, 011121112,\\ &&011221212, 011231324, 012221223, 012341224, 012342425 \}.
\end{eqnarray*} 

\noindent 3. $C_1:$ \begin{eqnarray*} \small   &&100000000, 110000000, 011000000, 001100000, 000110000,\\
 &&000010001, 000000111,
000001111, 000011101, 000110101,\\ &&001110011, 011110011, 111110101, 111111101,011111111, \\
&&001110112, 000120112, 000121112, 001111212, 011111212,\\
 &&111121112,
111220112, 012220112, 012221112, 111221212, \\
&&111121213, 011121223, 001221223,
001231213, 011231213,\\
 &&112221223, 122221223, 122231213, 112331213, 012331223, \\
&&012231224, 111231324, 111232324, 012232324, 012331324,\\
 &&112341224, 122341224,
123331324, 123332324, 122342324,\\ &&112341325, 012341335, 012342335, 112342425,
122342425,\\
 &&123342335, 123441335, 123451325, 123452325, 123442435,\\
 &&123342436,
122352436, 112452436, 013452436, 023452436.
 \end{eqnarray*}
 \noindent $\overset{\circ}{C} v_0 =v_1, \ldots, \overset{\bullet}{C} v_{59} =v_0 + \delta$.

\noindent 4. $C_2:$ \begin{gather*} \small     010000000, 111000000, 111100000, 011110000, 001110001,\\
 000110111, 000011112,
000011212, 000111212, 001121112,\\ 011220112, 112220112, 122221112, 122221212,
112221213,\\
 011231223, 001231224, 001231324, 011232324, 112232324,\\
 122331324,
123341224, 123441224, 123441324, 123342325,\\
 122342435, 112342436, 012352436,
012452436, 113452436.
\end{gather*}
 $\overset{\bullet}{C} v_0 =v_1, \ldots, \overset{\circ}{C} v_{29} =v_0 + \delta$.

\noindent 5. $C_3:$ \begin{gather*}\small
001000000, 011100000, 111110000, 111110001, 011110111,\\ 001111112, 000121212,
000121213, 001121223, 011221223,\\ 112231213, 122331213, 123331223, 123331224,
122341324,\\
 112342325, 012342435, 012342436, 112352436, 122452436.
\end{gather*}
 $\overset{\circ}{C} v_0 =v_1, \ldots, \overset{\bullet}{C} v_{19} =v_0 + \delta$.

\noindent 6. $C_4:$  \begin{eqnarray*} \small     &&000100000, 001110000, 011110001, 111110111, 111111112,\\ &&011121212, 001221213,
001231223, 011231224, 112231324,\\ &&122332324, 123342324, 123441325, 123451335,
123452336,\\
 &&123452536, 123453537, 123463547, 123562548, 124572548,\\ &&134673548,
235673648, 245673649, 245683659, 23578365(10),\\ &&13579375(10), 13579475(11),
23579486(11), 24579486(12), \\ &&2467(10)486(12).
\end{eqnarray*}
$\overset{\bullet}{C} v_0 =v_1, \ldots, \overset{\circ}{C} v_{29} =v_0 + 2 \delta$.

\noindent 7. $C_5:$ \begin{eqnarray*}\small    &&000010000, 000110001, 001110111, 011111112, 111121212,\\ &&111221213, 012231223,
012331224, 112341324, 122342325,\\ &&123342435, 123442436.
\end{eqnarray*}
$\overset{\circ}{C} v_0 =v_1, \ldots, \overset{\bullet}{C} v_{11} =v_0 +  \delta$.

\noindent 8. $C_6:$ \begin{gather*}\small
 000001000, 000001100, 000000101, 000010011, 000110011,\\
 001110101, 011111101,
111111111, 111110112, 011120112,\\ 001221112, 001221212,
 011121213, 111121223,
111221223,\\ 012231213, 012331213, 112331223, 122231224, 122231324,\\ 112332324, 012342324,
 012341325, 112341335, 122342335,\\ 123342425, 123442425,
 123452335,
123451336, 123451436.
 \end{gather*}

 $\overset{\bullet}{C} v_0 =v_1, \ldots, \overset{\circ}{C} v_{29} =v_0 +  \delta$.

\noindent 9. $C_7:$ \begin{eqnarray*} \small
 &&000000100, 000001101, 000011111, 000110112, 001120112,\\
 &&011221112, 112221212,
122221213, 122231223, 112331224,\\ &&012341324, 012342325, 112342435, 122342436, 123352436,\\
 &&123552436, 124562436, 134562437, 234562547, 234563548,\\
 &&134573648, 124673649, 124683659, 13468365(10),\\
 &&23568375(10), 24578475(10), 24679475(10), 24689375(11),\\
 &&2468(10)376(11), 2468(10)476(12).
\end{eqnarray*}
$\overset{\circ}{C} v_0 =v_1, \ldots, \overset{\bullet}{C} v_{29} =v_0 +  2 \delta$.

\noindent 10. $C_8:$  \begin{gather*}  000000010, 000000011, 000010101, 000111101, 001111111,\\
 011110112, 111120112,
111221112, 012221212, 012221213,\\ 111231223, 111231224, 012231324, 012332324
112342324,\\
 122341325, 123341335, 123442335, 123452425, 123452426.
\end{gather*}

$\overset{\circ}{C} v_0 =v_1, \ldots, \overset{\bullet}{C} v_{19} =v_0 +  \delta$.

\noindent 9. $C_0:$ \begin{gather*} \small    000000001, 000010111, 000111112, 001121212, 011221213,\\
 112231223, 122331224,
123341324, 123442325, 123452435.
\end{gather*}
 $\overset{\bullet}{C} v_0 =v_1, \ldots, \overset{\circ}{C} v_{9} =v_0 +  \delta$.

\noindent 10.

$$ \begin{array}{cc}
  \rm{A}_{1,k}= \qquad\qquad\qquad\qquad\qquad\qquad\qquad & \rm{A}_{2,k}= \qquad\qquad\qquad\qquad\qquad\qquad \\
   { D_{1, 65}\begin{cases}  C^{s-33} \overset{\bullet}{C} M_f  &
\text{if $k=2s$}\\
C^{s-33} M_f  & \text{if $k=2s-1$}
\end{cases}},  & { D_{2, 34}\begin{cases}  C^{s-17}  M_f  &
\text{if $k=2s$}\\
C^{s-17}\overset{\bullet}{C}  M_f  & \text{if $k=2s+1$}
\end{cases}};
 \\
\end{array}
 $$

$$ \begin{array}{cc}
  \rm{A}_{3,k}= \qquad\qquad\qquad\qquad\qquad\qquad\qquad & \rm{A}_{4,k}= \qquad\qquad\qquad\qquad\qquad\qquad \\
   { D_{3, 23}\begin{cases}  C^{s-12} \overset{\bullet}{C} M_f  &
\text{if $k=2s$}\\
C^{s-12}  M_f  & \text{if $k=2s-1$}
\end{cases}},  & { D_{4, 14}\begin{cases}  C^{s-7}  M_f  &
\text{if $k=2s$}\\
C^{s-7} \overset{\bullet}{C} M_f  & \text{if $k=2s+1$}
\end{cases}};
 \\
\end{array}
 $$

$$ \begin{array}{cc}
  \rm{A}_{5,k}= \qquad\qquad\qquad\qquad\qquad\qquad\qquad & \rm{A}_{6,k}= \qquad\qquad\qquad\qquad\qquad\qquad \\
  {D_{5, 13}\begin{cases}  C^{s-7} \overset{\bullet}{C} M_f  &
\text{if $k=2s$,}\\
C^{s-7}  M_f  & \text{if $k=2s-1$.}
\end{cases}},  & { D_{6, 28}\begin{cases}  C^{s-14}  M_f  &
\text{if $k=2s$}\\
C^{s-14}\overset{\bullet}{C}  M_f  & \text{if $k=2s+1$}
\end{cases}};
 \\
\end{array}
 $$

$$ \begin{array}{cc}
  \rm{A}_{7,k}= \qquad\qquad\qquad\qquad\qquad\qquad\qquad & \rm{A}_{8,k}= \qquad\qquad\qquad\qquad\qquad\qquad \\
   {D_{7, 17}\begin{cases}  C^{s-9} \overset{\bullet}{C} M_f  &
\text{if $k=2s$,}\\
C^{s-9}  M_f  & \text{if $k=2s-1$.}
\end{cases}},  & { D_{8, 19}\begin{cases}  C^{s-10} \overset{\bullet}{C} M_f  &
\text{if $k=2s$,}\\
C^{s-10}  M_f  & \text{if $k=2s-1$.}
\end{cases}};
 \\
\end{array}
 $$

$$\rm{A}_{0,k}=D_{0, 10}\begin{cases}  C^{s-5}  M_f  &
\text{if $k=2s$}\\
C^{s-5}  \overset{\bullet}{C}M_f  & \text{if $k=2s+1$}
\end{cases}
$$

$$ \begin{array}{cc}
  D_{1,65}= \qquad\qquad\qquad\qquad\qquad\qquad\qquad & D_{2,34} =\qquad\qquad\qquad\qquad\qquad\qquad \\
   \tiny{\left(%
\begin{array}{ccccccccc}
 -1 & 2 & -3 & 4 & -6 & 2 & -4 & -3 & 7 \\
 -2 & 4 & -6 & 8 & -11 & 4 & -7 & -5 & 12 \\
 -3 & 6 & -8 & 11 & -16 & 5 & -10 & -8 & 18 \\
 -3 & 7 & -10 & 14 & -21 & 7 & -14 & -10 & 24 \\
 -4 & 8 & -12 & 17 & -26 & 9 & -17 & -13 & 30 \\
 -2 & 3 & -5 & 7 & -10 & 3 & -7 & -5 & 12 \\
 -4 & 7 & -10 & 13 & -20 & 7 & -14 & -10 & 24 \\
 -2 & 5 & -8 & 10 & -15 & 5 & -10 & -8 & 18 \\
 -5 & 10 & -15 & 20 & -31 & 10 & -20 & -15 & 36 \\
\end{array}%
\right)},  & \tiny{\left(%
\begin{array}{ccccccccc
}
 -1 & 1 & -2 & 3 & -3 & 1 & -2 & -1 & 3 \\
 -1 & 2 & -4 & 5 & -6 & 2 & -4 & -3 & 7 \\
 -1 & 3 & -6 & 7 & -8 & 3 & -5 & -4 & 9 \\
 -2 & 4 & -7 & 8 & -10 & 4 & -7 & -5 & 12 \\
 -2 & 4 & -8 & 10 & -13 & 4 & -8 & -6 & 15 \\
 -1 & 2 & -3 & 4 & -5 & 1 & -3 & -3 & 6 \\
 -1 & 3 & -6 & 8 & -10 & 3 & -7 & -5 & 12 \\
 -1 & 2 & -5 & 6 & -7 & 2 & -5 & -4 & 9 \\
 -2 & 4 & -9 & 12 & -15 & 5 & -10 & -8 & 18 \\
\end{array}%
\right)};
 \\
\end{array}
 $$

$$ \begin{array}{cc}
  D_{3,23} = \qquad\qquad\qquad\qquad\qquad\qquad\qquad & D_{4, 14} = \qquad\qquad\qquad\qquad\qquad\qquad \\
   \tiny{ \left(%
\begin{array}{ccccccccc}
 -1 & 1 & -1 & 2 & -2 & 0 & -1 & -1 & 2 \\
 -1 & 2 & -3 & 4 & -4 & 1 & -2 & -2 & 4 \\
 -2 & 3 & -4 & 5 & -6 & 2 & -4 & -3 & 7 \\
 -3 & 4 & -5 & 6 & -7 & 3 & -5 & -3 & 8 \\
 -3 & 4 & -6 & 7 & -8 & 3 & -6 & -4 & 10 \\
 -1 & 2 & -2 & 2 & -3 & 1 & -2 & -2 & 4 \\
 -2 & 3 & -4 & 5 & -7 & 2 & -4 & -3 & 8 \\
 -1 & 2 & -3 & 4 & -5 & 2 & -4 & -2 & 6 \\
 -3 & 4 & -6 & 8 & -10 & 4 & -7 & -5 & 12 \\
\end{array}%
\right) },  & \tiny{ \left(%
\begin{array}{ccccccccc}
 0 & 1 & -1 & 1 & -1 & 0 & 0 & -1 & 1 \\
 -1 & 2 & -2 & 2 & -2 & 1 & -1 & -1 & 2 \\
 -1 & 2 & -3 & 3 & -3 & 1 & -2 & -2 & 4 \\
 -1 & 2 & -3 & 4 & -5 & 2 & -3 & -3 & 6 \\
 -1 & 2 & -3 & 4 & -6 & 2 & -3 & -3 & 7 \\
 0 & 0 & -1 & 2 & -2 & 1 & -1 & -1 & 2 \\
 0 & 1 & -2 & 3 & -4 & 2 & -3 & -2 & 5 \\
 -1 & 1 & -1 & 2 & -3 & 1 & -2 & -2 & 4 \\
 -1 & 2 & -3 & 4 & -6 & 3 & -4 & -4 & 8 \\
\end{array}%
\right) };
 \\
\end{array}
$$

$$ \begin{array}{cc}
  D_{5, 13} = \qquad\qquad\qquad\qquad\qquad\qquad\qquad & D_{6, 28} = \qquad\qquad\qquad\qquad\qquad\qquad \\
   \tiny{ \left(%
\begin{array}{ccccccccc}
 0 & 0 & 0 & 1 & -1 & 0 & -1 & 0 & 1 \\
 0 & 1 & -1 & 2 & -2 & 0 & -1 & -1 & 2 \\
 -1 & 2 & -2 & 3 & -3 & 1 & -2 & -1 & 3 \\
 -1 & 2 & -3 & 4 & -4 & 1 & -3 & -2 & 5 \\
 -1 & 2 & -3 & 5 & -6 & 2 & -4 & -3 & 7 \\
 0 & 0 & -1 & 2 & -2 & 1 & -1 & -1 & 2 \\
 0 & 1 & -2 & 3 & -4 & 2 & -3 & -2 & 5 \\
 -1 & 1 & -1 & 2 & -3 & 1 & -2 & -2 & 4 \\
 -1 & 2 & -3 & 5 & -7 & 2 & -4 & -3 & 8 \\
\end{array}%
\right) },  & \tiny{ \left(%
\begin{array}{ccccccccc}
-1 & 1 & -1 & 2 & -2 & 0 & -1 & -1 & 2 \\
 -1 & 2 & -3 & 4 & -4 & 1 & -2 & -2 & 4 \\
 -1 & 3 & -4 & 5 & -6 & 1 & -3 & -4 & 7 \\
 -2 & 4 & -5 & 6 & -8 & 2 & -5 & -5 & 10 \\
 -2 & 4 & -6 & 8 & -11 & 2 & -6 & -6 & 13 \\
 -1 & 2 & -3 & 4 & -5 & 1 & -3 & -3 & 6 \\
 -1 & 3 & -5 & 7 & -9 & 2 & -6 & -5 & 11 \\
 -1 & 2 & -4 & 5 & -6 & 1 & -4 & -4 & 8 \\
 -2 & 4 & -7 & 10 & -13 & 3 & -8 & -8 & 16 \\
\end{array}%
\right) };
 \\
\end{array}
 $$

 $$ \begin{array}{cc}
  D_{7, 17} = \qquad\qquad\qquad\qquad\qquad\qquad\qquad & D_{8, 19} = \qquad\qquad\qquad\qquad\qquad\qquad \\
   \tiny{ \left(%
\begin{array}{ccccccccc}
0 & 1 & -1 & 1 & -1 & 0 & 0 & -1 & 1 \\
 -1 & 2 & -2 & 2 & -2 & 1 & -1 & -1 & 2 \\
 -1 & 2 & -3 & 3 & -3 & 1 & -2 & -2 & 4 \\
 -1 & 2 & -3 & 4 & -5 & 2 & -3 & -3 & 6 \\
 -1 & 3 & -4 & 5 & -7 & 2 & -4 & -3 & 8 \\
 0 & 1 & -2 & 3 & -3 & 1 & -2 & -1 & 3 \\
 -1 & 3 & -4 & 5 & -6 & 2 & -4 & -3 & 7 \\
 -1 & 2 & -2 & 3 & -4 & 1 & -3 & -2 & 5 \\
 -1 & 4 & -5 & 6 & -8 & 3 & -6 & -4 & 10 \\
\end{array}%
\right) },  & \tiny{ \left(%
\begin{array}{ccccccccc}
-1 & 1 & -1 & 1 & -1 & 1 & -1 & 0 & 1 \\
 -1 & 1 & -2 & 2 & -2 & 1 & -2 & -1 & 3 \\
 -1 & 1 & -2 & 3 & -4 & 2 & -3 & -2 & 5 \\
 -1 & 2 & -3 & 4 & -6 & 2 & -4 & -2 & 7 \\
 -1 & 3 & -4 & 5 & -7 & 3 & -6 & -3 & 9 \\
 0 & 1 & -2 & 3 & -3 & 1 & -2 & -1 & 3 \\
 -1 & 2 & -3 & 5 & -6 & 2 & -5 & -2 & 7 \\
 -1 & 2 & -3 & 4 & -5 & 2 & -4 & -2 & 6 \\
 -2 & 4 & -5 & 7 & -9 & 3 & -7 & -4 & 11\\
\end{array}%
\right) };
 \\
\end{array}
 $$

 $$D_{0, 10} = \tiny{\left(%
\begin{array}{ccccccccc}
0 & 0 & 0 & 0 & -1 & 0 & 0 & 0 & 1 \\
 0 & 0 & 0 & 0 & -1 & 1 & -1 & -1 & 2 \\
 0 & 0 & 0 & 1 & -2 & 1 & -2 & -1 & 3 \\
 0 & 1 & -1 & 2 & -3 & 1 & -2 & -2 & 4 \\
 -1 & 2 & -2 & 3 & -4 & 2 & -3 & -2 & 5 \\
 -1 & 1 & -1 & 1 & -1 & 0 & -1 & -1 & 2 \\
 -1 & 1 & -2 & 3 & -3 & 1 & -2 & -2 & 4 \\
 0 & 1 & -2 & 2 & -2 & 1 & -2 & -1 & 3 \\
 -1 & 2 & -3 & 4 & -5 & 2 & -4 & -3 & 7\\
\end{array}%
\right)}
 $$
11.
 $\delta_{alg}=(1, 1, 1, 1, 1, 2, 2, 3, 6)$.

Hyperplane conditions $H_{\widetilde{E}_8}$:

$\delta + \alpha_2 > \gamma$, $\delta + \beta_1 > \gamma$, $\alpha_3 + \beta_1 > \gamma$, $3 \delta + \alpha_2 + \alpha_3 + \alpha_4 + \alpha_5 +   2 \beta_1 + 2 \beta_2 > 5 \gamma$,  $\delta + \alpha_1 + \alpha_5 + \beta_1 >
    2 \gamma$, $\alpha_1 + \beta_2 > \gamma$,  $\delta + \alpha_4 + \beta_2 > \gamma$,  $\delta + \alpha_2 + \alpha_4 + \beta_1 >
    2 \gamma$,  $\delta + \alpha_2 + \beta_1 + \beta_2 >     2 \gamma$,   $2 \delta + \alpha_2 + \alpha_3 + \beta_1 + \beta_2 >
    3 \gamma$,   $2 \delta + \alpha_1 + \alpha_4 + \beta_1 + \beta_2 >     3 \gamma$, $\delta + \alpha_1 + \alpha_3 + \alpha_4 + \beta_1 + \beta_2 >
    3 \gamma$, $\delta + \alpha_1 + \alpha_2 + \alpha_5 + \beta_1 + \beta_2 >
    3 \gamma$,  $\delta + \alpha_3 + \alpha_5 + \beta_1 + \beta_2 > 2 \gamma$,
  $2 \delta + \alpha_1 + \alpha_2 + \alpha_3 +\alpha_5 + \beta_1 + \beta_2 > 4 \gamma$, $2 \delta + \alpha_2 + \alpha_4 + \alpha_5 + \beta_1 + \beta_2 > 3 \gamma$,
  $2 \delta + \alpha_2 + \alpha_3 + \alpha_4 +
      2 \beta_1 + \beta_2 > 4 \gamma$, $2 \delta + \alpha_1 + \alpha_3 + \alpha_5 + 2 \beta_1 + \beta_2 > 4 \gamma$,
  $3 \delta + \alpha_1 + \alpha_2 + \alpha_4 + \alpha_5 +
      2 \beta_1 + \beta_2 > 5 \gamma$, $2 \delta + \alpha_1 + \alpha_2 + \alpha_3 + \alpha_4 + \alpha_5 + 2 \beta_1 + \beta_2 > 5 \gamma$,
  $2 \delta + \alpha_1 + \alpha_2 + \alpha_4 + \beta_1 +
      2 \beta_2 > 4 \gamma$,
  $3 \delta + \alpha_1 + \alpha_2 + \alpha_3 +\alpha_4 + \beta_1 + 2 \beta_2 > 5 \gamma$, $2 \delta + \alpha_1 + \alpha_3 + \alpha_4 +\alpha_5 + \beta_1 + 2 \beta_2 > 4 \gamma$,
  $2 \delta + \alpha_1 + \alpha_2 + \alpha_3 + \alpha_5 +
      2 \beta_1 + 2 \beta_2 > 5 \gamma$,  $\delta +\alpha_1 + \alpha_3 + \beta_2 >  2 \gamma$

\begin{center}
{\bf Acknowledgment}
\end{center}

The author thanks  Prof. Yu. Samoilenko and Prof. S. Kruglyak for many valuable
discussions and helpful comments. The paper was finished  at Chalmers
University of Technology, the author is grateful for the hospitality of the mathematics department and for the support of STINT. 
The author thanks the referee for recommending various
 improvements in exposition.


\end{document}